\documentclass[journal,twoside,web]{ieeecolor}
\pdfoutput=1
\usepackage{generic}

\makeatletter\let\NAT@parse\undefined\makeatother 
\usepackage{graphicx}
\usepackage{amsmath,amssymb}
\usepackage{mathtools}
\usepackage{cite}
\usepackage{xcolor}
\usepackage{booktabs,multirow}
\usepackage[flushleft]{threeparttable}

\usepackage{caption}
\usepackage[hidelinks]{hyperref}
\hypersetup{
	pdftitle={Schiller2026 - Performance guarantees for optimization-based state estimation using turnpike properties},
	pdfauthor={Julian D. Schiller, Matthias A. Mueller},
	colorlinks=true,
	linkcolor=cyan!60!black,
	citecolor=cyan!60!black,
	urlcolor=cyan!60!black,
}

\newtheorem{ass}{Assumption}
\newtheorem{defn}{Definition}
\newtheorem{thm}{Theorem}
\newtheorem{prop}{Proposition}
\newtheorem{lem}{Lemma}
\newtheorem{cor}{Corollary}
\newtheorem{rem}{Remark}

\newenvironment{pf}[1] {\noindent\hspace{2em}{\itshape {#1}: }}{\hspace*{\fill}~\QED\par\endtrivlist\unskip}

\usepackage{fancyhdr}
\def\doi{10.1109/TAC.2026.3755562}
\fancyfootoffset{3em}
\fancypagestyle{copyright}{
	\chead{
		\parbox[t]{0.82\headwidth}{%
			\centering\footnotesize
			This article has been accepted for publication in IEEE Transactions on Automatic Control. This is the author's version which has not been fully edited and content may change prior to final publication. Citation information: \href{http://dx.doi.org/\doi}{DOI \doi}
		}
	}
	\cfoot{%
		\raisebox{5mm}[0pt][0pt]{%
    		\parbox[t]{0.735\headwidth}{%
				\centering\footnotesize
				\textcopyright{}2026 IEEE. Personal use of this material is permitted.  Permission from IEEE must be obtained for all other uses, in any current or future media, including reprinting/republishing this material for advertising or promotional purposes, creating new collective works, for resale or redistribution to servers or lists, or reuse of any copyrighted component of this work in other works.
			}%
		}%
	}
}

\title{Performance guarantees for optimization-based state estimation using turnpike properties
\thanks{This work was supported by the Deutsche Forschungsgemeinschaft (DFG, German Research Foundation), Projects 426459964 and 499435839. \textit{(Corresponding author: Julian D. Schiller)}}
\thanks{Julian D. Schiller and Matthias A. Müller are with the Leibniz University Hannover, Institute of Automatic Control, 30167 Hannover, Germany. (e-mail: $\{$schiller,mueller$\}$@irt.uni-hannover.de).}
\thanks{Lars Grüne is with the University of Bayreuth, Chair of Applied Mathematics, 95447 Bayreuth, Germany. (e-mail: lars.gruene@uni-bayreuth.de).}}
\author{\parbox{\linewidth}{\centering Julian D. Schiller, Lars Grüne, and Matthias A. Müller}}

\begin{document}
	
\maketitle
\thispagestyle{copyright}
\pagestyle{empty}


\begin{abstract}	
	In this paper, we develop novel accuracy and performance guarantees for optimal state estimation of general nonlinear systems (in particular, moving horizon estimation, MHE). Our results rely on a turnpike property of the optimal state estimation problem, which essentially states that the omniscient infinite-horizon solution involving all past and future data serves as turnpike for the solutions of finite-horizon estimation problems involving a subset of the data. This leads to the surprising observation that MHE problems naturally exhibit a leaving arc, which may have a strong negative impact on the estimation accuracy. To address this, we propose a delayed MHE scheme, and we show that the resulting performance (both averaged and non-averaged) is approximately optimal and achieves bounded regret with respect to the infinite-horizon solution, with error terms that can be made arbitrarily small by an appropriate choice of the delay. In various simulation examples, we observe that already a very small delay in the MHE scheme is sufficient to significantly improve the overall estimation error by 20--25\% compared to standard MHE (without delay). This finding is of great importance for practical applications (especially for monitoring, fault detection, and parameter estimation) where a small delay in the estimation is rather irrelevant but may significantly improve the estimation results.
\end{abstract}

\section{Introduction}
Reconstructing the internal state trajectory of a dynamical system based on a batch of measured input-output data is an important problem of high practical relevance.
This can be accomplished, for example, by solving an optimization problem to find the best state and disturbance trajectories that minimize a suitably defined cost function depending on the measurement data.
If all available data is taken into account, this corresponds to the \emph{full information estimation} (FIE) problem.
However, if the data set or the underlying model is very large or only a limited amount of computation time or resources is available (as is the case with, e.g., online state estimation), the optimal solution to the FIE problem is usually difficult (or even impossible) to compute in practice.
For this reason, it is essential to find a reasonable approximation, which can be done, e.g., by means of a sequence of truncated optimal estimation problems, each of which uses only a limited time window of the full data set.
In the case of online state estimation, this corresponds to \emph{moving horizon estimation} (MHE), where in each discrete time step an optimal estimation problem is solved with a data set of fixed size.
MHE offers several advantages over classical state observers and filtering approaches due to its formulation as a constrained (deterministic) optimization problem. In particular, it is naturally applicable to nonlinear systems, allows for the incorporation of prior knowledge (e.g., on admissible state and disturbance trajectories), provides an intuitive tuning procedure, and enjoys strong theoretical properties, without being restricted to specific noise distributions.

\subsubsection*{Related work}
Current research in the field of nonlinear MHE is primarily concerned with stability and robustness guarantees, see, e.g., \cite[Ch.~4]{Rawlings2017} and \cite{Allan2021a,Knuefer2023,Schiller2023c,Hu2023,Alessandri2025}.
These works establish that, under a suitable detectability condition, the estimation error of MHE (i.e., the deviation between the estimated state and the real unknown system state) converges to a neighborhood of the origin, the size of which is determined by the magnitude of the unknown disturbance.
However, results on the actual \emph{performance} of MHE, and in particular on the approximation accuracy and performance loss (i.e., the \emph{regret}) with respect to a suitably challenging benchmark estimator, remain largely unavailable. This gap is primarily due to the intrinsic analytical complexity of nonlinear optimal estimation problems, for which general tools capable of yielding meaningful performance guarantees are not yet established.

For the special case of linear systems, regret-optimal filters for state estimation are designed in~\cite{Goel2023,Sabag2022} by minimizing the regret with respect to a clairvoyant (acausal) benchmark filter that has access to future measurements. This approach is extended in \cite{Brouillon2023}, where an exact solution to the minimal-regret observer is provided utilizing the system level synthesis framework; compare also \cite{Didier2022,Martin2024} for similar approaches in the context of regret-optimal control.
In \cite{Gharbi2020a}, an MHE scheme is proposed that provides regret guarantees with respect to an arbitrary comparative (e.g., the clairvoyant) observer. This is extended to nonlinear systems in~\cite{Gharbi2021}, but requires a restrictive convexity condition on the problem and disturbance- and noise-free data.

Whereas performance guarantees for state estimators are generally rather rare and usually restricted to linear systems, they often play an important role in nonlinear optimal control, especially when the overall goal is an economic one.
Corresponding results usually employ a \emph{turnpike} property of the underlying optimal control problem, cf.~\cite{McKenzie1986,Carlson1991}.
This property essentially implies that optimal trajectories stay close to an optimal equilibrium most of the time (or in general an optimal time-varying reference), which is regarded as the \emph{turnpike}.
Turnpike-related arguments are an important tool for assessing the closed-loop performance of model predictive controllers with general economic costs on finite and infinite horizons, cf. \cite{Gruene2016b,Faulwasser2022,Gruene2019}.
Necessary and sufficient conditions for the presence of the turnpike phenomenon in optimal control are discussed in, e.g.,~\cite{Damm2014,Gruene2016a,Trelat2023,Faulwasser2022a} and are usually based on dissipativity, controllability, and suitable optimality conditions.
\subsubsection*{Contributions}

In this paper, we establish novel accuracy and performance guarantees for nonlinear MHE, thereby addressing a fundamental gap in the literature. To this end, we develop systematic tools for the analysis of general nonlinear optimal state estimation problems based on turnpike arguments.
	
In particular, we characterize a novel turnpike property in this setting, whereby solutions to finite-horizon optimal estimation problems remain, for the majority of the horizon, close to the solution of the omniscient (acausal) infinite-horizon problem that incorporates all past and future data, referred to as the \emph{turnpike} (Section~\ref{sec:turnpike}).
This insight reveals that MHE optimization problems inherently exhibit a leaving arc at the end of the horizon, which can significantly degrade estimation accuracy.
To mitigate this phenomenon, we propose introducing an artificial delay in the MHE scheme, thereby realizing a fixed-lag moving horizon smoother (Section~\ref{sec:online}).
{We show that the resulting scheme is approximately optimal and achieves bounded regret with respect to the omniscient infinite-horizon solution, with error terms that can be made arbitrarily small by an appropriate choice of the delay.
Consequently, delayed MHE is capable of \emph{tracking the accuracy and performance of the omniscient infinite-horizon estimator} (Sections~\ref{sec:MHE:delay} and~\ref{sec:MHE_performance}).}
{Possible areas of application are those in which a small delay between receiving the measurement and estimating the state is tolerable (which is generally a requirement of smoothing methods), e.g., navigation, target tracking, audio signal processing, biological processes, and parameter learning, compare also \cite[Sec.~1.1]{Saerkkae2013}.}

As additional contribution, we propose a novel turnpike prior for MHE formulations with prior weighting, \mbox{which---in} contrast to the standard filtering or smoothing \mbox{prior---can} be shown to converge to a neighborhood around the turnpike that can be made arbitrarily small by design \mbox{(Section~\ref{sec:MHE:prior})}.
Furthermore, we consider the case of batch (i.e., offline) state estimation (Section~\ref{sec:offline}) and show that good performance of a state estimator directly implies small estimation errors (Section~\ref{sec:accuracy}).
We illustrate our theoretical findings with several simulation examples from the literature (Sections~\ref{sec:example}), including a continuously stirred tank reactor and a highly nonlinear quadrotor model with 12 states (the code is publicly available, see~\cite{Schiller2025_code}). In all examples, we observe the turnpike phenomenon being present and find that even a delay of very few steps in the MHE scheme improves the overall estimation error by 20\,--25\% compared to standard MHE (without delay).

Compared to the preliminary conference version~\cite{Schiller2024b}, in this paper we focus on the case of MHE for online state estimation. Specifically, we provide novel results on the accuracy and performance of MHE with respect to the omniscient infinite-horizon estimator incorporating all past and future data, establish bounded regret with respect to this benchmark, and consider the practically relevant case of MHE formulations with an additional prior weighting.
Moreover, we link estimator performance to its accuracy with respect to the true unknown system trajectory, give further technical details, and provide more extensive simulation examples. Additional results on turnpikes in optimal state estimation can also be found in~\cite{Schiller2025a}.

\subsubsection*{Notation}
We denote the set of integers by $\mathbb{Z}$, the set of natural numbers (including zero) by $\mathbb{N}$, and the set of \emph{even} natural numbers by $\mathbb{N}_{\mathrm{e}}$.
When we use interval notation, we always refer to integer intervals.
The notation $a\pm b$ for $a,b\in\mathbb{N}$ refers to the symmetric interval $[a-b,a+b]$.
We denote a sequence $\{x_j\}_{j=a}^b$ of elements $x_j\in\mathbb{R}^n$, $j\in[a,b]$ with $x_{a:b}$.
The Euclidean norm of a vector $x\in\mathbb{R}^n$ is denoted by $|x|$, and the weighted norm with respect to a symmetric positive definite matrix $Q\succ0$ by $|x|_Q=\sqrt{x^\top Q x}$. The maximum eigenvalue of $Q$ is denoted by $\lambda_{\max}(Q)$.
The identity matrix of size $n\times n$ is denoted by $I_n$.
Finally, we recall that a function $\alpha:\mathbb{R}_{\geq 0}\rightarrow\mathbb{R}_{\geq 0}$ is of class $\mathcal{K}$ if it is continuous, strictly increasing, and satisfies $\alpha(0)=0$.
By $\mathcal{L}$, we refer to the class of functions $\theta:\mathbb{R}_{\geq 0}\rightarrow \mathbb{R}_{\geq 0}$ that are continuous, non-increasing, and satisfy $\lim_{s\rightarrow\infty}\theta(s)=0$, and by $\mathcal{KL}$ to the class of functions $\beta:\mathbb{R}_{\geq 0}\times\mathbb{R}_{\geq 0}\rightarrow\mathbb{R}_{\geq 0}$ with $\beta(\cdot,s)\in\mathcal{K}$ and $\beta(r,\cdot)\in\mathcal{L}$ for any fixed $s\in\mathbb{R}_{\geq 0}$ and $r\in\mathbb{R}_{\geq 0}$, respectively.

\section{Problem Setup}\label{sec:problem}

\subsection{System description}\label{sec:system_description}
We consider nonlinear uncertain discrete-time systems:
\begin{subequations}\label{eq:sys}
	\begin{align}
		{x}_{t+1} &= f(x_t,u_t,w_t),\label{eq:sys_1}\\
		y_t &= h(x_t,u_t) + v_t\label{eq:sys_2}
	\end{align}
\end{subequations}
with discrete time $t\in\mathbb{N}$, state $x_t\in\mathbb{R}^n$, (known) control input $u_t\in\mathbb{R}^m$, (unknown) process disturbance $w_t\in\mathbb{R}^q$, (unknown) measurement noise $v_t\in\mathbb{R}^p$, and noisy output measurement $y_t\in\mathbb{R}^p$.
The functions $f:\mathbb{R}^n\times\mathbb{R}^m\times\mathbb{R}^q\rightarrow \mathbb{R}^n$ and $h:\mathbb{R}^n\times\mathbb{R}^m\rightarrow \mathbb{R}^p$ define the system dynamics and output equation.

We assume that trajectories of the system~\eqref{eq:sys} satisfy
\begin{equation}\label{eq:con_Z}
	(x_t,u_t,w_t,v_t) \in \mathcal{X}\times\mathcal{U}\times\mathcal{W}\times\mathcal{V},\  {t\in\mathbb{N}}
\end{equation}
for some known sets $\mathcal{X}\subseteq\mathbb{R}^n$, $\mathcal{U}\subseteq\mathbb{R}^m$, $\mathcal{W}\subseteq\mathbb{R}^q$, $\mathcal{V}\subseteq\mathbb{R}^p$, and furthermore, that $(x,u,w)\in\mathcal{X}\times\mathcal{U}\times\mathcal{W}\Rightarrow f(x,u,w)\in\mathcal{X}$.
Such knowledge typically arises from the physical nature of the system (e.g., non-negativity of partial pressure or molar concentration), the incorporation of which can significantly improve the estimation results, cf.~\cite[Sec.~4.4]{Rawlings2017} and \cite{Schiller2023c} and compare the examples in Section~\ref{sec:example}. If no such information is available, one can simply choose $\mathcal{X}=\mathbb{R}^n$, $\mathcal{U}=\mathbb{R}^m$, $\mathcal{W}=\mathbb{R}^q$, $\mathcal{V}=\mathbb{R}^p$, thus trivially satisfying the above~assumptions.

\subsection{The optimal state estimation problem}
For ease of notation, we define the input-output data tuple $d_t := (u_t,y_t)$ obtained from the system~\eqref{eq:sys} at time $t\in\mathbb{N}$. Consider a given batch of input-output data $d_{0:N}$ of length $N+1$ for some $N\in\mathbb{N}$.
We aim to compute the state and disturbance sequences $\hat{x}_{0:N}$ and $\hat{w}_{0:N-1}$ that are optimal in the sense that they minimize a cost function involving the full data set $d_{0:N}$ (see its definition in \eqref{eq:MHE_cost} below).
In particular, we consider the following optimal state estimation problem
\begin{subequations}\label{eq:MHE}
	\begin{align}\label{eq:MHE_0}
		P_N(d_{0:N}):\hspace{5.0ex} &\min_{\text{\footnotesize$\begin{matrix}
					\hat{x}_{0:N}\\ \hat{w}_{0:N-1}
				\end{matrix}$}}\  J_{N}(\hat{x}_{0:N},\hat{w}_{0:N-1}; d_{0:N}) \\ 
		\text{s.t. }\
		&\hat{x}_{j+1}=f(\hat{x}_{j},u_j,\hat{w}_{j}), \ j\in[0,N-1], \label{eq:MHE_f} \\
		&\hat{x}_{j}\in\mathcal{X}, \ j\in[0,N], \label{eq:MHE_x} \\
		&\hat{w}_{j}\in\mathcal{W},\ j\in[0,N-1], \label{eq:MHE_w}\\
		& y_j-h(\hat{x}_{j},u_j)\in\mathcal{V}, \ j\in[0,N].  \label{eq:MHE_y}
	\end{align}
\end{subequations}
For convenience, we define the combined sequence $\hat{z}_{0:N}$ as $\hat{z}_j:=(\hat{x}_j,\hat{w}_j)$ for $j\in[0,N-1]$ and $\hat{z}_N := (\hat{x}_N,0)$. The constraints~\eqref{eq:MHE_f}--\eqref{eq:MHE_y} enforce the prior knowledge about the system model, the domain of the true trajectories, and the disturbances and noise (note that feasibility is always guaranteed due to the fact that the true trajectory satisfies~\eqref{eq:con_Z} by assumption).
In~\eqref{eq:MHE_0}, we consider the cost function
\begin{align}
	J_N(\hat{x}_{0:N},\hat{w}_{0:N-1} ;d_{0:N})  = \sum_{j=0}^{N-1}l(\hat{x}_{j},\hat{w}_{j};d_j) + g(\hat{x}_{N};d_N)\label{eq:MHE_cost}
\end{align}
with continuous stage cost $l:\mathcal{X}\times\mathcal{W}\times\mathcal{U}\times\mathbb{R}^p\rightarrow\mathbb{R}_{\geq0}$ and terminal cost $g:\mathcal{X}\times\mathcal{U}\times\mathbb{R}^p\rightarrow\mathbb{R}_{\geq0}$.
This is a generalization of classical designs used for state estimation, where $l$ and $g$ are positive definite in the disturbance input $\hat{w}$ and the fitting error $y-h(\hat{x},u)$, cf.,~e.g.,~\cite[Ch.~4]{Rawlings2017}; it particularly includes the practically relevant case of quadratic costs
\begin{align}
	l(x,w;d) &= |w|^2_Q + |y-h(x,u)|^2_R,\label{eq:stage_cost}\\
	g(x;d) &= |y-h(x,u)|^2_G, \label{eq:term_cost}
\end{align}
where $Q,R,G$ are symmetric positive definite matrices.
However, our results also hold for more general functions $l$ and $g$, allowing the objective in~\eqref{eq:MHE_cost} to be tailored to the specific problem at hand. Designs involving an additional prior weighting in~\eqref{eq:MHE_cost} (as usual in MHE) are considered in Section~\ref{sec:MHE:prior}.

In the state estimation context, a cost function with terminal cost $g$ is usually referred to as the \textit{filtering form} of the state estimation problem.
To simplify the analysis, most recent works on MHE theory consider a cost function without terminal cost (i.e., $g=0$), which corresponds to the \textit{prediction form} of the state estimation problem~\cite{Schiller2023c,Knuefer2023,Allan2021a,Hu2023,Alessandri2025}, compare~\cite[Ch.~4]{Rawlings2017}.
Overall, the formulation~\eqref{eq:MHE}--\eqref{eq:MHE_cost} (together with the extension in Section~\ref{sec:MHE:prior}) covers all optimal estimation problems underlying most modern MHE schemes used in theory and practice.

The optimization problem $P_N$~\eqref{eq:MHE} is a parametric nonlinear program, the solution of which solely depends on the (input-output) data provided, i.e., the sequence $d_{0:N}$.
We characterize solutions to $P_N$ using the generic solution mapping $\zeta_N$:
\begin{equation}
	\hat{z}_{j}^* := \zeta_N(j,d_{0:N}),\ j\in[0,N], \label{eq:z_N}
\end{equation}
which yields the value function $V_N(d_{0:N})=J_N(\hat{z}_{0:N}^{*};d_{0:N})$.
While $\zeta_N$ returns a pair of optimal state and disturbance input, we use $\zeta_N^x$ to refer only to the optimal state of this pair, so that $\hat{x}^*_j = \zeta_N^x(j,d_{0:N})$ for all $j\in[0,N]$.

\begin{rem}[Existence of solutions to $P_N$]\label{rem:existence}
	We assume that, whenever we employ the solution mapping $\zeta_N$, the corresponding solution exists. Under continuity of $f$ and $h$, this can generally be guaranteed by choosing the stage cost $l$ and terminal cost~$g$ such that the cost function $J_N$ is radially unbounded in the (condensed) decision variables (which requires positive definiteness of $l$ and $g$ and observability\footnote{
		Observability is required here because the cost function $J_N$ does not contain an additional prior weighting; the case of MHE with prior weighting, which does not require observability, is considered in Section~\ref{sec:MHE:prior} below.
	} of the system with a corresponding horizon length~$N$, compare~\cite[Sec.~4.3.1]{Rawlings2017}), or under compactness of the sets $\mathcal{X}$ and $\mathcal{W}$, see \cite[Prop.~A.7]{Rawlings2017}.
\end{rem}

\subsection{Benchmark solution}

We are interested in how the solution in~\eqref{eq:z_N} compares to a certain (challenging) benchmark problem.
For this purpose, we interpret the measured data sequence $d_{0:N}$ as a segment of an infinite data sequence $d_{-\infty:\infty}$ that contains all past and future data that could possibly be generated by the system~\eqref{eq:sys} on $\mathbb{Z}$.
Then, we consider the \emph{omniscient infinite-horizon estimator}, that is, the solution of the (acausal) optimal state estimation problem
\begin{equation}\label{eq:IHE}
	 P_\infty(d_{-\infty:\infty}) \ : \ \min\limits_{\hat{z}_{-\infty:\infty}} \sum_{j=-\infty}^{\infty}l(\hat{z}_{j};d_j) \ \ \text{s.t.}\ \  \eqref{eq:MHE_f}\text{--}\eqref{eq:MHE_y}, \ j\in\mathbb{Z},
\end{equation}
where $\hat{z}_j=(\hat{x}_j,\hat{w}_j)$, $j\in\mathbb{Z}$. We denote the solution to $P_\infty(d_{-\infty:\infty})$ by $z^\infty:=z^\infty_{-\infty:\infty}$ with $z^\infty_j=(x^\infty_j,w^\infty_j)$, $j\in\mathbb{Z}$, where we assume that $z^\infty$ exists and is unique for any possible $d_{-\infty:\infty}$, compare Remark~\ref{rem:existence}. Extending our theory to the case where $z^\infty$ is non-unique is an interesting topic for future research, compare the discussion at the end of Section~\ref{sec:turnpike_sensitivity}.

\begin{rem}[Infinite-horizon optimality]
	In general, the objective function in~\eqref{eq:IHE} does not attain a finite value. In this case, we consider~$z^\infty$ to be optimal in the sense of the concept of \emph{overtaking optimality}, compare Definitions 1 and 2 in~\cite{Gruene2019}.
\end{rem}

{In linear settings \cite{Sabag2022,Brouillon2023}, a common benchmark for observers estimating the state $x_t$, $t\in\mathbb{N}$ is the clairvoyant (acausal) observer relying on the infinite data sequence $d_{0:\infty}$, where in particular future data $d_{t+1:\infty}$ is only fictitious (as it depends on future disturbances and noise) and may or may not actually be measured at a later point in time (e.g., if the experiment is terminated).}
We extend this concept by assuming that our benchmark, the omniscient infinite-horizon estimator, has both access to the future and a perfect memory.

In the following, we investigate how the finite-horizon solution $\hat{z}^*_{0:N}$ behaves compared to the infinite-horizon benchmark solution $z^\infty$ on the interval $[0,N]$, quantifying how much we loose by \emph{not} being able to compute the latter in practice.

\section{Turnpike in optimal state estimation problems}\label{sec:turnpike}

Optimal state estimation problems (such as~\eqref{eq:MHE}) can be interpreted as optimal control problems using the disturbance~$\hat{w}$ as the control input, cf.~\cite[Sec.~4.2.3]{Rawlings2017}, \cite[Sec.~4]{Allan2020a}.
In particular, a cost function~\eqref{eq:MHE_cost} that is positive definite in the estimated disturbance and the fitting error (as, e.g., in \eqref{eq:stage_cost} and \eqref{eq:term_cost}) can be regarded as an economic output tracking cost, penalizing deviations from the ideal reference $(w^\mathrm{r}_j,y^\mathrm{r}_j) = (0,y_j)$, $j\in[0,N]$.
This reference, however, is generally unreachable, as it is usually impossible to attain zero cost $V_N(d_{0:N})=0$, except for the special case where $y_{0:N}$ corresponds to an output sequence of~\eqref{eq:sys} under zero disturbances $w_{0:N-1}\equiv0$, $v_{0:N}\equiv0$.
For unreachable references, on the other hand, it is known that the corresponding optimal control problem exhibits the turnpike property with respect to the \emph{best reachable reference}~\cite{Koehler2019}, which suggests that a similar phenomenon can also be expected in optimal state estimation problems.

In the following Section~\ref{sec:mot_example}, we provide a simple example that supports this intuition. Then, we mathematically formalize and discuss this property in Section~\ref{sec:turnpike_sensitivity}.

\subsection{Motivating example}\label{sec:mot_example}
Suppose that some output data $y_{0:T}$ for $T=30$ is measured from the system $x_{t+1} = x_t + w_t$, $x_0=1$, $y_t = x_t + v_t$, where $w_t = v_t = 1$ for all $t\in\mathbb{Z}$.
We consider the finite-horizon optimal estimation problem $P_N(d_{0:N})$ with quadratic costs~\eqref{eq:stage_cost} and~\eqref{eq:term_cost} using \mbox{$Q=R=G=1$} and different values of $N$.
For comparison reasons, we also consider the infinite-horizon problem $P_\infty$, which we approximate by simulating the system and computing the solution on some extended interval $[-T_e,T+T_e]$, where we choose $T_e$ such that the corresponding solution does not change on $[0,T]$ if $T_e$ is further increased.

\begin{figure}
	\hspace{0.4ex}
	\includegraphics{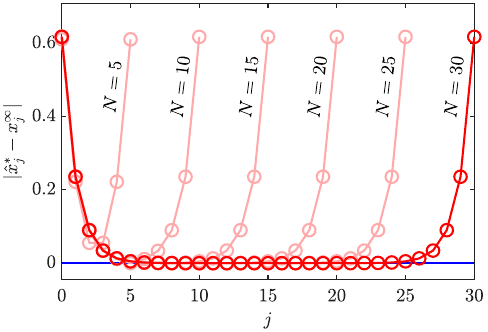}
	\caption{Difference between the solution of the infinite-horizon problem $P_\infty$ and solutions of the finite-horizon problem $P_N$ for different values of $N$.}
	\label{fig:OE}
\end{figure}

Figure~\ref{fig:OE} shows the difference between the infinite-horizon solution $x^\infty_j$ and the solution of the finite-horizon problem $\hat{x}^*_j$ for $j\in[0,N]$. 
{Here, we find that the sequence $\hat{x}^*_{0:N}$ generally consists of three pieces: first, a transient at the beginning where $\hat{x}^*_j$ converges to ${x}^\infty_j$; second, a large middle phase where $\hat{x}^*_j$ stays near ${x}^\infty_j$; third, a transient at the end where $\hat{x}^*_j$ diverges from ${x}^\infty_j$. Figure~\ref{fig:OE} indicates that these transients are independent of the horizon length $N$. Overall, it can be concluded that, when $N$ is large, the solution $\hat{x}^*_j$ remains close to ${x}^\infty_j$ for most of $j\in[0,N]$. This is in fact the key characteristic of any turnpike definition appearing in the literature,~e.g.,~\cite{Gruene2016b,Faulwasser2022,Gruene2019,Gruene2016a,Trelat2023}.}
{From this perspective, we can refer to the solution ${x}^\infty_j$ as the \emph{turnpike}, the left transient as the \emph{approaching arc} (or \emph{entry arc}), and the right transient as the \emph{leaving arc}, cf.~\cite{Faulwasser2022}.}

\subsection{Turnpike characterization}\label{sec:turnpike_sensitivity}

There are in fact multiple ways to formalize the turnpike behavior observed in the motivating example in Section~\ref{sec:mot_example}.
A characterization popular in the context of (receding horizon) optimal control involves a bound (independent of the horizon length) on the number of elements of the sequence $\hat{z}^*_{0:N}$ that lie outside an \mbox{$\epsilon$-neighborhood} of the turnpike, cf.~\cite{Gruene2016a,Gruene2018}.
{This \emph{measure turnpike} property is particularly suitable for use in economic model predictive control (cf., e.g.,~\cite{Faulwasser2018}), and also has the decisive advantage that the corresponding sufficient condition that is usually invoked in this context (namely, \emph{strict dissipativity}) is a global concept.}
{In fact, we can easily adapt the standard result that ``strict dissipativity
implies turnpike'' (e.g., \cite[Th.~1]{Gruene2018}) to the more general
setting of optimal state estimation, see \cite[Th.~6.1]{Schiller2025a} for details.}
Unfortunately, the resulting turnpike property is rather weak in the sense that it is not possible to infer \emph{which} elements of the solution $\hat{z}^*_{0:N}$ are actually close to the turnpike, and which are not.
Such additional information (which is crucially required in the context of state estimation, as will be clear in Sections~\ref{sec:online} and~\ref{sec:offline}) is provided by exponential (or polynomial) turnpike characterizations that involve an explicit time-dependent bound on the difference of optimal trajectories and the turnpike, cf.~e.g., \cite{Damm2014,Trelat2015,Trelat2023,Han2022}.
To cover arbitrary decay rates, we propose the following unified turnpike property involving general $\mathcal{KL}$-functions.

\begin{defn}[Turnpike for optimal state estimation]\label{def:turnpike}
	\, The optimal estimation problem $P_N$ exhibits the $\mathcal{KL}$-turnpike property with respect to the infinite-horizon solution $z^\infty$ if there exists $\beta\in\mathcal{KL}$ such that $\hat{z}^*_{j} = \zeta_N(j,d_{0:N})$ satisfies
	\begin{equation}\label{eq:turnpike}
		|\hat{z}_{j}^*-z_{j}^\infty|\leq \beta(|\hat{x}_{0}^*-x_{0}^\infty|,j) + \beta(|\hat{x}_{N}^*-x_{N}^\infty|,N-j)
	\end{equation}
	for all $j\in[0,N]$, $N\in\mathbb{N}$, and all possible data $d_{-\infty:\infty}$.
\end{defn}

The two terms on the right-hand side in~\eqref{eq:turnpike} capture the occurrence of approaching and leaving arcs. Intuitively, these reflect the fact that the finite-horizon problem $P_N$ considers only a segment of the full data set that underlies the infinite-horizon problem $P_\infty$, neglecting both the fictitious historical data $d_{-\infty:-1}$ and the future data $d_{N+1:\infty}$.

{For the special case of linear systems and quadratic cost functions as in~\eqref{eq:stage_cost} and~\eqref{eq:term_cost}, in~\cite[Sec.~6.2]{Schiller2025a} we show that the turnpike property from Definition~\ref{def:turnpike} is naturally present under observability and controllability (the latter one with respect to the disturbance input) using standard arguments from optimal control and Riccati theory. In fact, using Proposition 6.2 and Theorem~6.2 in \cite{Schiller2025a}, the $\mathcal{KL}$-function $\beta$ in \eqref{eq:turnpike} specializes to $\beta(r,s) = Kr\lambda^s$ for some $K>0$ and $\lambda\in(0,1)$, constituting an exponential turnpike property.}
{Hence, the conditions are intuitive, not particularly restrictive in the context of state estimation, and can be easily checked.}

\begin{rem}[Motivating example revisited]\label{rem:motivation_TP}
	{As the system considered in Section~\ref{sec:mot_example} is linear and trivially observable and controllable, the corresponding optimal state estimation problem exhibits the turnpike property from Definition~\ref{def:turnpike} (with exponential decay).
	In line with our theory, Figure~\ref{fig:OE} shows clear approaching and leaving arcs for any choice of~$N$.}
\end{rem}

For nonlinear systems, however, deriving sufficient conditions for turnpike behavior in the spirit of Definition~\ref{def:turnpike} usually requires a detailed and rather technical analysis.
For example, one may establish~\eqref{eq:turnpike} by considering linearizations of the extremal equations as in \cite{Trelat2015} or by exploiting the specific graph structure of the underlying nonlinear program and use the decaying sensitivity results from \cite{Na2020,Shin2022}.
For the latter, it has been shown that exponentially decaying sensitivity is present under standard regularity and optimality conditions of the problem (such as a uniform second order sufficient condition and uniform boundedness of the Lagrangian Hessian, cf.~\cite{Shin2022}), which are satisfied under uniform observability and controllability conditions, see~\cite{Shin2021} and compare also \cite{Gruene2020}.
However, these approaches suffer from local validity and may only hold in a possibly small neighborhood of the turnpike.

{Here, it is worth noting that establishing \emph{global} (exponential) turnpike properties of optimal control problems is generally difficult, particularly due to the potential non-uniqueness of the turnpike, see, e.g.,~\cite{Pighin2022} for further details on this topic and compare our discussion below~\eqref{eq:IHE}.
Corresponding works rely on combining assumptions of global nature (such as strict dissipativity) with assumptions of local nature that involve the linearizations at the turnpike, cf.~\cite{Damm2014,Trelat2023}. Another line of research deals with Lipschitz nonlinear systems and uses bootstrap arguments to establish global exponential turnpike behavior, thereby avoiding global dissipativity assumptions and the use of linearizations, cf.~\cite{EsteveYaguee2022}.
Generalizing these results to our case is an interesting topic for future research.}

In the following, we use the turnpike property from Definition~\ref{def:turnpike} to assess the performance of MHE and its regret with respect to the infinite-horizon benchmark $z^\infty$.

\section{Optimization-based online state estimation}\label{sec:online}
\begin{figure}[t]
	\centering
	\includegraphics[width=\columnwidth]{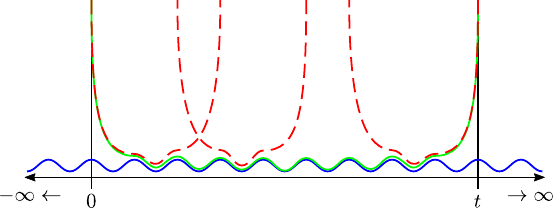}
	\caption{Sketch of the infinite-horizon solution $x^\infty$ (blue), the current FIE solution $\zeta_t(j,d_{0:t})$ (green), and solutions of the finite-horizon estimation problem $\zeta_N(j,d_{\tau-N:\tau})$ for different values of $\tau$ (red).}
	\label{fig:turnpike_sol}
\end{figure}

In online state estimation, one is generally interested in obtaining, at each time instant $t\in\mathbb{N}$, an accurate estimate of the current true (unknown) state $x_t$.
A natural approach is to simply solve the optimal state estimation problem in~\eqref{eq:MHE} based on all available historical data (by setting $N=t$). This corresponds to the case of FIE, which can be formalized using the solution mapping $\zeta_N^x$ defined below~\eqref{eq:z_N} as follows:
\begin{equation}
	\hat{x}_{t}^{\mathrm{fie}} = \zeta_t^x(t,d_{0:t}), \ t\in\mathbb{N}. \label{eq:FIE_seq_x}
\end{equation}
However, repeatedly solving $P_t(d_{0:t})$ for the current FIE solution $\hat{x}^{\mathrm{fie}}_t$ is generally infeasible in practice since the problem size continuously grows with time.
Instead, MHE considers the truncated optimal estimation problem $P_N(d_{t-N:t})$ using the most recent data $d_{t-N:t}$ only, where the horizon length $N\in\mathbb{N}$ is fixed.
More precisely, the MHE estimate at time $t\in\mathbb{N}$ can be formulated as
\begin{equation}
	\hat{x}_{t}^{\mathrm{mhe}} =
	\begin{cases}
		\zeta_N^x(N,d_{t-N:t}), & t\in[N,\infty)\\
		\zeta_t^x(t,d_{0:t}), & t \in[0,N-1].
	\end{cases}\label{eq:MHE_seq_x}
\end{equation}
MHE constitutes a well-established method for state estimation that is increasingly applied in practice. However, assuming that the underlying optimal estimation problem exhibits the turnpike property in the sense of Definition~\ref{def:turnpike}, we know that both the FIE sequence defined by~\eqref{eq:FIE_seq_x} and the MHE sequence defined by~\eqref{eq:MHE_seq_x} consist of point estimates of finite-horizon problems that lie on the leaving arc, see Figure \ref{fig:turnpike_sol}; hence, \emph{FIE and MHE might produce estimates that are actually far away from the turnpike}.

To mitigate this, we propose the following modifications:
\begin{itemize}
	\item Reduce the influence of the leaving arc by introducing an artificial delay in the estimation, cf. Sections~\ref{sec:MHE:delay}--\ref{sec:MHE_performance}.
	\item Reduce the influence of the approaching arc by using a turnpike-based prior weighting, cf.~Section~\ref{sec:MHE:prior}.
\end{itemize}

The subsequent analyses rely on the assumption that the MHE optimization problems satisfy a turnpike property.

\begin{ass}\label{ass:turnpike_inf}
	The finite-horizon optimal state estimation problem $P_N$ in~\eqref{eq:MHE} has the turnpike property (Definition~\ref{def:turnpike}).
\end{ass}

Verifying Assumption~\ref{ass:turnpike_inf} is generally challenging in non-convex optimal state estimation problems, see Section~\ref{sec:turnpike_sensitivity}. In practice, this may be assessed empirically by solving the estimation problem for varying lengths $N$ and examining whether the resulting solutions exhibit the turnpike phenomenon; see Section~\ref{sec:example}, where it can be observed even for strongly nonlinear systems. Importantly, the proposed modifications do not require verification of Assumption~1 for implementation.

In the remainder of this section, we consider $N\in\mathbb{N}_{\mathrm{e}}$.

\subsection{A delay improves the estimation results} \label{sec:MHE:delay}
To reduce the influence of naturally appearing leaving arcs in MHE problems~\eqref{eq:MHE_seq_x}, we introduce a delayed MHE scheme ($\delta$MHE). For a fixed delay $\delta \in [0,N/2]$, $N\in\mathbb{N}_\mathrm{e}$, we define
\begin{equation}
	\hat{x}_{t-\delta}^{\mathrm{\delta mhe}} =
	\begin{cases}
		\zeta_N^x(N-\delta,d_{t-N:t}), & t\in[N,\infty)\\
		\zeta_t^x(t-\delta,d_{0:t}), & t \in [\delta,N-1].
	\end{cases}\label{eq:MHE_delay_x}
\end{equation}
Note that for the special case $\delta=0$, $\delta$MHE reduces to standard MHE~\eqref{eq:MHE_seq_x}.
Hence, the parameter $\delta$ is an additional degree of freedom that constitutes a trade-off between delaying the estimation results and reducing the influence of the leaving~arc.

\begin{rem}[$\delta$FIE]
	By replacing $N$ with $t$ and using only the second case in \eqref{eq:MHE_delay_x}, we can similarly define a delayed FIE scheme---thus, all of the following results and implications derived for $\delta$MHE directly carry over to $\delta$FIE.
\end{rem}

\begin{rem}[Smoothing form of MHE]\label{rem:smoothing}
	Considering a fixed delay in the estimation scheme to improve the results is actually quite common in signal processing and filtering theory and refers to \emph{fixed-lag smoothing algorithms} \cite[Ch.~5]{Crassidis2011}, \cite{Saerkkae2013}.
	Early results for linear systems can be found in, e.g., \cite{Moore1973}; more recent works address, e.g., robustness against model errors \cite{Yi2023} or extensions to certain classes of nonlinear systems \cite{Rehman2015}.
	For linear systems with Gaussian noise, fixed-lag receding-horizon smoothers are proposed in \cite{Simon2013,Kim2013,Kwon2022}.
	In this context, the proposed $\delta$MHE scheme can also be interpreted as the (fixed-lag) smoothing form of MHE; compared to the smoothing literature, however, we consider general nonlinear systems under arbitrary process and measurement noise and provide performance and regret guarantees with respect to the infinite-horizon optimal solution using turnpike~arguments.
\end{rem}

The following result shows that by a suitable choice of the delay $\delta$, the state estimates~\eqref{eq:MHE_delay_x} can be driven arbitrarily close to the omniscient infinite-horizon benchmark estimator $z^\infty$.

\begin{prop}\label{prop:dMHE}
	Let Assumption~\ref{ass:turnpike_inf} hold and $\mathcal{X}$ be compact. Then, there exists $\sigma\in\mathcal{L}$ such that the estimated state sequence of $\delta$MHE in~\eqref{eq:MHE_delay_x} satisfies
	\begin{equation}\label{eq:delta_mhe_turnpike}
		|\hat{x}^{\delta \mathrm{mhe}}_{j}-x^\infty_{j}|\leq \sigma(\delta), \ j\in[\delta,t-\delta]
	\end{equation}
	for all $t\in[\delta,\infty)$, $\delta\in[0,N/2]$, $N\in\mathbb{N}_\mathrm{e}$, and all $d_{-\infty:\infty}$.
\end{prop}
\begin{proof}
	By compactness of $\mathcal{X}$, there exists $C>0$ such that $|x_1-x_2|\leq C$ for all $x_1,x_2\in\mathcal{X}$.
	The rest of this proof follows directly from the definition of $\delta$MHE~\eqref{eq:MHE_delay_x} and Assumption~\ref{ass:turnpike_inf}.
	Specifically, for all $j\in[\delta,N-\delta-1]$, we have
	\begin{align*}
		|\hat{x}^{\delta\mathrm{mhe}}_j - x^\infty_j| &\leq 	|\zeta_{j+\delta}(j,d_{0:j+\delta})-z^\infty_j|\\
		&\leq \beta(|\zeta_{j+\delta}^x(0,d_{0:j+\delta})-x^\infty_0|,j) \\
		&\quad +  \beta(|\zeta_{j+\delta}^x(j+\delta,d_{0:j+\delta})-x^\infty_{j+\delta}|,\delta)\\
		&\leq 2\beta(C,\delta).
	\end{align*}	
	For $j\in[N-\delta,t-\delta]$, on the other hand, it follows that
	\begin{align*}
		|\hat{x}^{\delta\mathrm{mhe}}_j - x^\infty_j|&\leq 	|\zeta_N(N-\delta,d_{j+\delta-N:j+\delta})-z^\infty_j|\\
		&\leq \beta(|\zeta_N^x(0,d_{j+\delta-N:j+\delta})-x^\infty_{j+\delta-N},N-\delta) \\
		&\quad +  \beta(|\zeta_N^x(N,d_{j+\delta-N:j+\delta})-x^\infty_{j+\delta}|,\delta)\\
		&\leq 2\beta(C,\delta).
	\end{align*}	
	Combining both cases for $j\in[\delta,t-\delta]$ and defining $\sigma(s) := 2\beta(C,s)$, $s\geq0$ yields \eqref{eq:delta_mhe_turnpike}. Noting that $\sigma\in\mathcal{L}$ establishes the statement and hence finishes this proof.
\end{proof}

Proposition~\ref{prop:dMHE} provides an estimate on the accuracy of $\delta$MHE in the sense of how close the sequence $\hat{x}^{\delta\mathrm{mhe}}_{\delta:t-\delta}$ is to the benchmark ${x}^\infty_{\delta:t-\delta}$. Since $\sigma$ in~\eqref{eq:delta_mhe_turnpike} is of class $\mathcal{L}$, this difference can be made arbitrarily small by increasing the delay $\delta$ (as well as $N$ if needed). This is in line with intuition, as increasing $\delta$ results in the state estimates in~\eqref{eq:MHE_delay_x} being closer to the turnpike. Note that under an exponential turnpike property, $\sigma$ is of \emph{exponential} form, see Remark~\ref{rem:exp_turnpike} below for further details.

\subsection{Performance estimates}\label{sec:MHE_performance}

{In this section, for simplicity we restrict our analysis to dynamics~\eqref{eq:sys_1} that involve additive disturbances (cf.~Remark~\ref{rem:nonlinear_w}):}
\begin{equation}
	f(x,u,w) = f_\mathrm{a}(x,u) + w \label{eq:dynamics}
\end{equation}
with $w\in\mathcal{W}=\mathbb{R}^n$.
Furthermore, we impose a Lipschitz condition on the nonlinear functions $f_\mathrm{a}$ and $h$ on $\mathcal{X}$.

\begin{ass}\label{ass:Lipschitz}
	The functions $f_\mathrm{a}$ and $h$ are Lipschitz in $x\in\mathcal{X}$, i.e., there exist constants $L_f,L_h>0$ such that
	\begin{align}
		|f_\mathrm{a}(x_1,u)-f_\mathrm{a}(x_2,u)| &\leq L_f|x_1-x_2|, \label{eq:Lipschitz_f}\\
		|h(x_1,u)-h(x_2,u)| &\leq L_h|x_1-x_2| \label{eq:Lipschitz_h}
	\end{align}
	for all $x_1,x_2\in\mathcal{X}$ uniformly for all $u\in\mathcal{U}$.
\end{ass}

Note that Assumption~\ref{ass:Lipschitz} is not restrictive in practice under compactness of $\mathcal{X}$ (as considered in Proposition~\ref{prop:dMHE}).

The dynamics~\eqref{eq:dynamics} ensure one-step controllability with respect to the disturbance input $w$; consequently, the sequence $\hat{x}^{\delta\mathrm{mhe}}_{\delta:t-\delta}$ forms a state trajectory of system~\eqref{eq:sys}, driven by the disturbance input $\hat{w}_j^{\delta\mathrm{mhe}} = \hat{x}^{\delta\mathrm{mhe}}_{j+1}- f_\mathrm{a}(\hat{x}^{\delta\mathrm{mhe}}_{j},u_j)$, $j\in[\delta,t-\delta-1]$.
For the sake of conciseness, we define $\hat{z}^{\delta\mathrm{mhe}}_j := (\hat{x}^{\delta\mathrm{mhe}}_j,\hat{w}^{\delta\mathrm{mhe}}_j)$, $j\in[\delta,t-\delta-1]$ and  $\hat{z}^{\delta\mathrm{mhe}}_{t-\delta} :=(\hat{x}^{\delta\mathrm{mhe}}_{t-\delta},0)$.

{We now specify our measure of performance, relying on the stage cost $l$ from~\eqref{eq:MHE_cost}.
Let $t_1,t_2\in\mathbb{N}$ denote some interval of interest $[t_1,t_2]$.}
For a given sequence $\hat{z}_{t_1:t_2}$ with $\hat{z}_j = (\hat{x}_j,\hat{w}_j)$ satisfying the system dynamics $\hat{x}_{j+1} = f(\hat{x}_j,u_j,\hat{w}_j)$ for $j\in[t_1,t_2]$, we consider the performance criterion
\begin{equation}\label{eq:perf_criterion}
	J_{[t_1,t_2]}(\hat{z}_{t_1:t_2}):=\sum_{j=t_1}^{t_2-1}l(\hat{z}_j;d_j).
\end{equation}

\begin{rem}[{Nonlinear disturbances}]\label{rem:nonlinear_w}
	{The proof of Theorem~\ref{thm:perf_MHE} below (and hence all results in this section) can be generalized to nonlinear disturbances as in~\eqref{eq:sys_1} (i.e., omitting conditions~\eqref{eq:dynamics} and \eqref{eq:Lipschitz_f}) by defining the performance criterion~\eqref{eq:perf_criterion} for general sequences $\hat{x}_{t_1:t_2}$ that do \emph{not} obey the system dynamics~\eqref{eq:sys_1}, e.g., $J_{[t_1,t_2]}(\hat{x}_{t_1:t_2})=\sum_{j=t_1}^{t_2-1} |\hat{x}_{j+1}-\hat{x}_j|^2_Q + |y_j-h(\hat{x}_j,u_j)|^2_R$ for some $Q,R\succ0$.
	{Yet, this way we lose the direct link to the optimization problem~\eqref{eq:MHE}.}}
\end{rem}

{We provide the following performance estimate for $\delta$MHE.}

\begin{thm}[Performance of $\delta$MHE]\label{thm:perf_MHE}
	Consider the system dynamics \eqref{eq:dynamics} and the quadratic stage cost in~\eqref{eq:stage_cost}.
	Let Assumptions \ref{ass:turnpike_inf} and \ref{ass:Lipschitz} be satisfied and $\mathcal{X}$ be compact.
	Then, there exists $\bar{\sigma}\in\mathcal{L}$ such that for any choice of $\epsilon>0$, $\delta$MHE for some arbitrary delay $\delta\in[0,N/2]$ and horizon length $N\in\mathbb{N}_\mathrm{e}$ satisfies the following performance estimate:
	\begin{equation}
		J_{[t_1,t_2]}(\hat{z}^\mathrm{\delta mhe}_{t_1:t_2}) \leq (1+\epsilon)J_{[t_1,t_2]}({z}^{\infty}_{t_1:t_2}) + \frac{1+\epsilon}{\epsilon} t_\Delta\bar{\sigma}(\delta) \label{eq:thm_performance}
	\end{equation}
	for all $t_1,t_2\in[\delta,t-\delta]$, all $t\in[\delta,\infty)$, and all possible $d_{-\infty,\infty}$, where $t_\Delta=t_2-t_1$.
\end{thm}

Before proving Theorem~\ref{thm:perf_MHE}, we want to highlight some key properties of the performance estimate in~\eqref{eq:thm_performance}.

\begin{rem}[Performance of $\delta$MHE] \label{rem:imterpretation} \phantom{c}
	\begin{enumerate}
		\item The performance of $\delta$MHE (the sequence $\hat{z}^{\delta\mathrm{mhe}}_{t_1:t_2}$) is approximately optimal on the interval $[t_1,t_2]$ with error terms that depend on $\epsilon$, $\delta$, and the interval length $t_\Delta$.
		
		\item The performance estimate in \eqref{eq:thm_performance} directly implies a bound on the regret (i.e., the performance loss) of $\delta$MHE with respect to the benchmark $z^\infty$, see~Corollary~\ref{cor:regret}.
		
		\item {In case of an exponential turnpike property, the performance of $\delta$MHE converges \emph{exponentially} to the performance of the infinite-horizon estimator as $\delta$ increases. For further details, see Remark~\ref{rem:exp_turnpike} below.}
		
		\item The performance estimate~\eqref{eq:thm_performance} grows linearly with the size of the performance interval $t_\Delta$ and tends to infinity if $t_\Delta$ approaches infinity. This property is to be expected (due to the fact that the turnpike is never exactly reached) and conceptually similar to (non-averaged) performance results in economic model predictive control, see, e.g., \cite[Sec.~5]{Faulwasser2018}, \cite{Gruene2016b}. To make meaningful performance estimates in case $t_\Delta\rightarrow\infty$, we analyze the averaged performance in Corollary~\ref{cor:performance}.
		
		\item Theorem~\ref{thm:perf_MHE} is stated for the practically relevant case of quadratic stage costs as in~\eqref{eq:stage_cost} for ease of presentation, but can easily be extended to more general cost functions that fulfill a weak triangle inequality, compare also~\cite{Schiller2024b}.
	\end{enumerate}
\end{rem}

\begin{pf}{Proof of Theorem~\ref{thm:perf_MHE}}
	Using the definitions from~\eqref{eq:MHE_cost}--\eqref{eq:term_cost}, the performance criterion evaluated for $\delta$MHE reads
	\begin{equation}
		J_{[t_1,t_2]}(\hat{z}^{\delta\mathrm{mhe}}_{t_1:t_2})= \sum_{j=t_1}^{t_2-1}|\hat{w}^{\delta\mathrm{mhe}}_{j}|_Q^2 + |y_j - h(\hat{x}^{\delta\mathrm{mhe}}_{j},u_j)|_R^2. \label{eq:proof_JT}
	\end{equation}
	From the definition of $\hat{x}^{\delta\mathrm{mhe}}_j$ in~\eqref{eq:MHE_delay_x} and the fact that $f_\mathrm{a}(x_{j}^\infty,u_j)+ w^\infty_j - x_{j+1}^\infty=0$ using~\eqref{eq:dynamics}, the square root of the input cost satisfies
	\begin{align}
		&|\hat{w}^{\delta\mathrm{mhe}}_{j}|_Q\nonumber\\
		& = |\hat{x}^{\delta\mathrm{mhe}}_{j+1}- x_{j+1}^\infty + f_\mathrm{a}(x_{j}^\infty,u_j) - f_\mathrm{a}(\hat{x}^{\delta\mathrm{mhe}}_{j},u_j) + w^\infty_j|_Q\nonumber\\
		& \leq \sqrt{\lambda_{\max}(Q)}(|\hat{x}^{\delta\mathrm{mhe}}_{j+1}- x_{j+1}^\infty| + L_f|x_{j}^\infty-\hat{x}^{\delta\mathrm{mhe}}_{j}|) + |w^\infty_j|_Q \label{eq:proof_wc}
	\end{align}
	for all $j\in[t_1,t_2-1]$, where in the last step we have used the triangle inequality and Assumption~\ref{ass:Lipschitz}.
	From Proposition~\ref{prop:dMHE} and the fact that $t_1,t_2\in[\delta,t-\delta]$, it follows that
	\begin{equation}
		|\hat{x}^{\delta\mathrm{mhe}}_j - x^\infty_j| \leq |\hat{z}^{\delta\mathrm{mhe}}_j - z^\infty_j| \leq \sigma(\delta), \ j\in[t_1,t_2],\label{eq:proof_perf_x}
	\end{equation}
	where $\sigma\in\mathcal{L}$.
	By squaring both sides of~\eqref{eq:proof_wc} and using~\eqref{eq:proof_perf_x} together with the fact that for any $\epsilon>0$, it holds that $(a + b)^2\leq (1+\epsilon)a^2 + \frac{1+\epsilon}{\epsilon}b^2$ for all $a,b\geq0$, we obtain
	\begin{equation}
		|\hat{w}^{\delta\mathrm{mhe}}_{j}|^2_Q \leq (1+\epsilon)|w_j^\infty|^2_Q + \frac{1+\epsilon}{\epsilon}(1+L_f)^2\lambda_{\max}(Q)\sigma(\delta)^2 \label{eq:proof_Q}
	\end{equation}
	for each $j\in[t_1,t_2-1]$.
	A similar reasoning for the fitting error (where we add $h(x_j^\infty,u_j) - h(x_j^\infty,u_j) = 0$) yields
	\begin{equation*}
		|y_j - h(\hat{x}^{\delta\mathrm{mhe}}_{j},u_j)|_R\leq |y_j - h(x_j^\infty,u_j)|_R + L_h\sqrt{\lambda_{\max}(R)}\sigma(\delta)
	\end{equation*}
	for all $j\in[t_1,t_2-1]$. By squaring both sides and using the same argument that allowed us to obtain~\eqref{eq:proof_Q}, we get
	\begin{align}
		|y_j - h(\hat{x}^{\delta\mathrm{mhe}}_{j},u_j)|_R^2 
		&\leq (1+\epsilon)|y_j - h(x_j^\infty,u_j)|_R^2 \nonumber \\
		&\quad + \frac{1+\epsilon}{\epsilon}L_h^2\lambda_{\max}(R)\sigma(\delta)^2 \label{eq:proof_R}
	\end{align}
	for all $j\in[t_1,t_2-1]$.
	The performance criterion~\eqref{eq:proof_JT} together with~\eqref{eq:proof_Q}, \eqref{eq:proof_R}, and the definition of $\bar{\sigma}(s):=\big((1+L_f)^2\lambda_{\max}(Q)+L_h^2\lambda_{\max}(R) \big) \sigma(s)^2,\, s\geq0$ (where we note that $\bar{\sigma}\in\mathcal{L}$) establishes the statement of this theorem.
\end{pf}

To derive a bound on the regret and the asymptotic averaged performance of $\delta$MHE, we first show that $J_{[t_1,t_2]}({z}^{\infty}_{t_1:t_2})$ grows at maximum linearly in $t_\Delta$.

\begin{lem}\label{lem:VT}
	Let $\mathcal{W}$ and $\mathcal{V}$ be compact. There exists $A>0$ such that $J_{[t_1,t_2]}({z}^{\infty}_{t_1:t_2}) \leq  A(t_2-t_1)$ for any possible $d_{-\infty:\infty}$.
\end{lem}
\begin{proof}
	Due to $\mathcal{W}$ and $\mathcal{V}$ being compact, there exist $C_{\mathrm{Q}},C_{\mathrm{R}}>0$ such that $|\hat{w}|^2_Q\leq C_{\mathrm{Q}}$ and $|y-h(\hat{x},u)|_R^2\leq C_{\mathrm{R}}$ for all $(\hat{x},\hat{w},u,y)\in\mathcal{X}\times\mathcal{W}\times\mathcal{U}\times\mathbb{R}^p$ such that $y-h(\hat{x},u)\in\mathcal{V}$. Hence, the claim holds with $A:=C_{\mathrm{Q}}+C_{\mathrm{R}}$.
\end{proof}

The following corollary from Theorem~\ref{thm:perf_MHE} establishes bounded regret of $\delta$MHE with respect to the benchmark $z^\infty$.

\begin{cor}[Bounded regret]\label{cor:regret}
	Let the conditions of Theorem~\ref{thm:perf_MHE} be satisfied.
	Assume that $\mathcal{W}$ and $\mathcal{V}$ are compact.
	Then, the regret of $\delta$MHE can be bounded by
	\begin{equation*}
		J_{[t_1,t_2]}({z}^{\delta\mathrm{mhe}}_{t_1:t_2}) - J_{[t_1,t_2]}(z^\infty_{t_1:t_2}) \leq t_{\Delta}\left(\epsilon A + \frac{1+\epsilon}{\epsilon}\bar{\sigma}(\delta)\right)
	\end{equation*}
	for all $t_1,t_2\in[\delta,t-\delta]$, all $t\in[\delta,\infty)$, and all possible $d_{-\infty,\infty}$, where $\epsilon>0$ and $\bar{\sigma}\in\mathcal{L}$ are from Theorem~\ref{thm:perf_MHE}, $A>0$ is from Lemma~\ref{lem:VT}, and $t_\Delta=t_2-t_1$.
\end{cor}

Note that the regret bound provided by Corollary~\ref{cor:regret} is linear in the interval length $t_{\Delta}$, where the slope $C(\epsilon,\delta):=\left(\epsilon A + \frac{1+\epsilon}{\epsilon}\bar{\sigma}(\delta)\right)$ can be rendered arbitrarily small by suitable choices of $\epsilon$ and $\delta$. Again, linear dependency on $t_\Delta$ is to be expected as the turnpike is never exactly reached, cf.~Remark~\ref{rem:imterpretation}.

The following result establishes an estimate on the averaged performance of $\delta$MHE for the asymptotic case when $t_{\Delta}\rightarrow\infty$.

\begin{cor}[Averaged performance]\label{cor:performance}
	Let the conditions of Theorem~\ref{thm:perf_MHE} be satisfied.
	Assume that $\mathcal{W}$ and $\mathcal{V}$ are compact.
	Then, $\delta$MHE satisfies the averaged performance estimate
	\begin{equation*}
		\limsup\limits_{t_\Delta\rightarrow\infty}\frac{1}{t_\Delta}(J_{[t_1,t_2]}({z}^{\delta\mathrm{mhe}}_{t_1:t_2}) - J_{[t_1,t_2]}({z}^\infty_{t_1:t_2})) \leq \epsilon A + \frac{1+\epsilon}{\epsilon}\bar{\sigma}(\delta)
	\end{equation*}
	for all possible data $d_{-\infty,\infty}$, where $\epsilon>0$ and $\bar{\sigma}\in\mathcal{L}$ are from Theorem~\ref{thm:perf_MHE}, $A>0$ is from Lemma~\ref{lem:VT}, and $t_\Delta=t_2-t_1$.
\end{cor}

Corollary~\ref{cor:performance} implies that the averaged performance of the $\delta$MHE is finite (due to the fact that $J_{[t_1,t_2]}({z}^\infty_{t_1:t_2})$ can be bounded as $J_{[t_1,t_2]}({z}^\infty_{t_1:t_2})/t_\Delta\leq A$ by Lemma~\ref{lem:VT}) and approximately optimal with respect to the infinite-horizon benchmark estimator, with error terms that can be made arbitrarily small by suitable choices of $\epsilon$ and $\delta$.
If we choose $\delta$ as an appropriate function of $\epsilon$, we can even fully recover the benchmark performance in the limit $\delta\rightarrow\infty$, as implied by the next~corollary.

\begin{cor}[Latency vs. performance trade-off]\label{cor:regret_2}
	{\phantom{x} Let the conditions of Theorem~\ref{thm:perf_MHE} be satisfied.
	Assume that $\mathcal{W}$ and $\mathcal{V}$ are compact. There exists a strictly decreasing function $\delta(\epsilon)$ with $\delta:\mathbb{R}_{>0}\rightarrow\mathbb{R}_{>0}$ and $\lim_{\epsilon\rightarrow0}\delta(\epsilon)=\infty$ together with some $\alpha\in\mathcal{K}$ such that
	\begin{equation*}
			\limsup\limits_{
					\text{\scriptsize$\begin{matrix} t_\Delta\rightarrow\infty\\ \end{matrix}$}
				} \frac{1}{t_\Delta}\left( J_{[t_1,t_2]}({z}^{\delta\mathrm{mhe}}_{t_1:t_2}) - J_{[t_1,t_2]}(z^\infty_{t_1:t_2})\right) \leq \alpha(\epsilon)
	\end{equation*}
	for any $\epsilon>0$ and all possible data $d_{-\infty,\infty}$, where $t_\Delta=t_2-t_1$.}
\end{cor}

\begin{proof}
	{Let $g(s):=\bar{\sigma}(s) + 1/s$, $s>0$, where $\bar{\sigma}\in\mathcal{L}$ is from Theorem~\ref{thm:perf_MHE}.
	Note that $g:\mathbb{R}_{>0}\rightarrow\mathbb{R}_{>0}$ is one-to-one; hence, its inverse $g^{-1}:\mathbb{R}_{>0}\rightarrow\mathbb{R}_{>0}$ exists.
	Define $\delta:\mathbb{R}_{>0}\rightarrow\mathbb{R}_{>0}$ as $\delta(s) = g^{-1}(s^2)$, $s>0$.
	As $g$ is strictly decreasing, $g^{-1}$ and $\delta$ are as well.
	Moreover, because $1/s$ is strictly positive for all $s>0$, it follows that $\bar{\sigma}(\delta(s)) < g(\delta(s)) = g(g^{-1}(s^2)) = s^2$. Using this property in the bound provided by Corollary~\ref{cor:regret} and taking the limit yields
	\begin{equation*}
		\limsup\limits_{t_\Delta\rightarrow\infty}\frac{1}{t_\Delta}(J_{[t_1,t_2]}({z}^{\delta\mathrm{mhe}}_{t_1:t_2}) - J_{[t_1,t_2]}({z}^\infty_{t_1:t_2})) \leq \epsilon^2 + (1+A)\epsilon.
	\end{equation*}
	Using $\alpha(s):=s^2 + (1+A)s$, $s>0$ concludes this proof.}
\end{proof}

Corollary~\ref{cor:regret_2} characterizes the trade-off between latency and performance, i.e., between a desired upper bound on the averaged regret of $\delta$MHE and the delay $\delta$ required for this. For the averaged regret approaching zero, we would need that $\delta\rightarrow\infty$, while for a larger bound smaller values of $\delta$ suffice.

Overall, Proposition~\ref{prop:dMHE}, Theorem~\ref{thm:perf_MHE}, and Corollaries~\ref{cor:regret}--\ref{cor:regret_2} imply that $\delta$MHE is able to track the solution and the performance of the omniscient infinite-horizon estimator.
Increasing the delay $\delta$ reduces the influence of the leaving arc and improves the performance estimate; the best performance is achieved for $\delta=N/2$, which, on the other hand, introduces a potentially large delay (depending on $N$).

\begin{rem}[{Performance with exponential turnpike}]\label{rem:exp_turnpike}
	{In the practically relevant case of exponential turnpike behavior, Assumption~\ref{ass:turnpike_inf} holds with $\beta(r,t)=Kr\lambda^{t}$ in Definition~\ref{def:turnpike} for some $K>0$ and $\lambda\in(0,1)$. Then, the $\mathcal{L}$-functions $\sigma$ used in Proposition~\ref{prop:dMHE} and $\bar{\sigma}$ used in Theorem~\ref{thm:perf_MHE} and Corollaries~\ref{cor:regret}--\ref{cor:regret_2} are decreasing exponentially---more precisely, $\sigma(s)=K_1\lambda^s$ and $\bar{\sigma}(s) = K_2\lambda^{2s},s\geq0$ for some $K_1,K_2>0$.}
	{Consequently, already small delays are expected to significantly reduce the influence of the leaving arc and hence improve the estimation results compared to standard MHE (without delay), which is empirically validated by our simulations in Section~\ref{sec:example}.}
\end{rem}

We conclude this section by noting that while it is possible in model predictive control to design suitable terminal ingredients that yield finite non-averaged performance for $t\rightarrow\infty$ (cf.~\cite{Grune2015}), this does not seem to be possible here, as it would imply that we have certain information about future data in order to exactly reach and stay on the solution of the acausal infinite-horizon estimation problem.

\subsection{MHE with prior weighting}\label{sec:MHE:prior}

It is known that MHE schemes with a cost function as in~\eqref{eq:MHE_cost} might require relatively large estimation horizons to achieve satisfactory estimation results, cf.~\cite[Sec.~4.3.2]{Rawlings2017}.
In order to reduce the required horizon length and enable faster computations, MHE formulations that leverage an additional prior weighting are therefore usually preferred in practice.
The prior weighting can generally be seen as additional regularization of the cost function, ensuring that the initial state $\hat{x}_0$ of an estimated sequence $\hat{x}_{0:N}$ stays in a meaningful region.
In view of our turnpike results, a well-chosen prior weighting hence reduces the influence of the approaching arc and ensures that solutions of truncated problems can reach the turnpike in fewer steps.

The prior weighting is usually parameterized by a given prior estimate $\bar{x}$ and a (possibly time-varying) function $p_t(x,\bar{x})$ that is positive definite and uniformly bounded in the difference $|x-\bar{x}|$.
Throughout the following, we use the definition $N_t := \min\{t,N\}$, which is convenient here as it avoids additional case distinctions.
The (time-varying) MHE cost function with prior weighting can then be formulated as
\begin{equation}\label{eq:J_p}
	\begin{split}
		&J_{N_t}^{\mathrm{p}}(\hat{x}_{0:N_t},\hat{w}_{0:N_t-1};d_{t-N_t:t},t)\\
		&:= p_{t-N_t}(\hat{x}_0,\bar{x}_{t-N_t})  + J_{N_t}(\hat{x}_{0:N_t},\hat{w}_{0:N_t-1};d_{t-N_t:t}),
	\end{split}
\end{equation}
with the current decision variables $\hat{x}_{0:N_t}$ and $\hat{w}_{0:N_t-1}$, and where $J_{N_t}$ is from~\eqref{eq:MHE_cost} with $N$ replaced by $N_t$.
At any time $t\in\mathbb{N}$, the current MHE problem to solve is given by the optimal estimation problem~\eqref{eq:MHE} with $N$ replaced by $N_t$ and the cost function $J_{N_t}^\mathrm{p}$ from~\eqref{eq:J_p}, which we denote by $P_{N_t}^\mathrm{p}(d_{t-N_t:t},\bar{x}_{t-N_t},t)$. The corresponding solution (which exists under mild conditions, compare Remark~\ref{rem:existence}) is described by the sequence $\tilde{z}^*_{t-N_t:t}$, where $\tilde{z}_{j}^* = (\tilde{x}_{j}^*,\tilde{w}_{j}^*)$, $j\in[t-N_t,t-1]$ and $\tilde{z}_{t}^* = (\tilde{x}_{t}^*,0)$.
The prior $\bar{x}_{t-N_t}$ is typically chosen in terms of a past solution of the problem $P_{N_t}^\mathrm{p}$, which introduces a coupling between the MHE problems.
For easier reference, it is therefore convenient to introduce double indices, where, e.g., $\tilde{z}_{j|t}^*$, $j\in[t-N_t,t]$ refers to the element $\tilde{z}_j^*$ of the solution of $P_{N_t}^\mathrm{p}(\cdot,t)$ computed at time $t\in\mathbb{N}$.
We first consider the standard MHE case and set the current estimate $\hat{x}^{\mathrm{mhe,p}}_t$ to the last state of the sequence $\tilde{z}_{t-N_t:t|t}^*$, i.e., $\hat{x}^{\mathrm{mhe,p}}_t := \tilde{x}_{t|t}^*$, $t\in\mathbb{N}$.
In Remark~\ref{rem:perf} below, we again consider a delay to reduce the influence of the naturally appearing leaving arc.

A popular choice for the prior weighting is
\begin{equation}\label{eq:prior_cost}
	p_t(x,\bar{x})=|x-\bar{x}|_{W_{t}}^2,
\end{equation}
where $W_t$ is a constant or time-varying positive definite weighting matrix that might be updated using, e.g., covariance update laws from nonlinear Kalman filtering, cf., \cite{Rao2003,Qu2009,Baumgaertner2020,Kuehl2011}.
Given some initial guess $\bar{x}_0\in\mathcal{X}$, there are two common choices for updating the prior estimate $\bar{x}_{t-N_t}$ at time $t\in\mathbb{N}$: first, the \emph{filtering prior}
\begin{equation}\label{eq:prior_filtering}
	\bar{x}_{t-N_t} =
	\begin{cases}
		\tilde{x}_{t-N|t-N}^*, & t\in[N,\infty)\\
		\bar{x}_{0}, & t\in[0,N-1],
	\end{cases}
\end{equation}
which corresponds to the state estimate computed $N$ steps in the past, i.e., the last element of the solution to the problem {$P_{N}^\mathrm{p}(\cdot,t-N)$}; second, the \emph{smoothing prior}
\begin{equation}\label{eq:prior_smoothing}
	\bar{x}_{t-N_t} = 
	\begin{cases}
		\tilde{x}_{t-N_t|t-1}^*, & t\in[1,\infty)\\
		\bar{x}_{0}, & t=0,
	\end{cases}
\end{equation}
which refers to the second element of the solution to $P^\mathrm{p}_{N_t}(\cdot,t-1)$ computed at the previous time step $t-1$, cf.~\cite[Sec.~4.3.2]{Rawlings2017}. MHE with a smoothing prior can generally recover faster from poor initial estimates, whereas MHE with a filtering prior essentially comprises measurements from two time horizons and may therefore be advantageous in the long term.
Note that the filtering prior is considered in most of the recent literature on nonlinear MHE theory, as this allows to derive a contraction of the estimation error over the estimation horizon and thus establish robust stability \cite{Allan2021a,Knuefer2023,Schiller2023c,Hu2023,Alessandri2025,Rawlings2017}.
In contrast, the smoothing prior is used in, e.g., \cite{Alessandri2008,Alessandri2017}, and often serves in practice-oriented works as a linearization point for computing an improved prior estimate using Kalman filter updates, compare, e.g., \cite{Kuehl2011,Baumgaertner2020} and see also the review article \cite{Elsheikh2021}.

However, from our turnpike analysis, we know that both of these two options may be unsuitable if we are interested in approximating the infinite-horizon optimal performance; the smoothing prior corresponds to an element of a finite-horizon solution on the approaching arc, and the filtering prior to an element of the solution on the leaving arc. Hence, we propose the following \emph{turnpike prior}:
\begin{equation}\label{eq:prior_turnpike}
	\bar{x}_{t-N_t} =
	\begin{cases}
		\tilde{x}_{t-N_t|t-N/2}^*, & t\in[N/2,\infty)\\
		\bar{x}_{0}, & t\in[0,N/2-1],
	\end{cases}
\end{equation}
which corresponds to the middle element of the solution of $P_{N}(\cdot,t-N/2)$ computed at time $t-N/2$. This avoids the influence of the approaching and leaving arcs (for $t\geq N$).
In fact, we can show that the prior estimate $\bar{x}_{t-N_t}$ converges to a neighborhood of the turnpike $x^\infty$ under the following modified (exponential) turnpike property of the MHE problem~$P_{N_t}^\mathrm{p}$.

\begin{ass}[Turnpike for MHE with prior weighting]\label{ass:turnpike_inf_prior}
	There exist constants $K>0$ and $\lambda\in(0,1)$ such that for all $N\in\mathbb{N}$, the solutions of the finite-horizon problem $P_N^\mathrm{p}(d_{\tau:\tau+N},\bar{x},\tau+N)$ and the infinite-horizon problem $P_\infty(d_{-\infty:\infty})$ satisfy
	\begin{equation*}
		|\tilde{z}_{\tau+j}^*-z^\infty_{\tau+j}|\leq K\left(|\tilde{x}_{\tau}^*-x^\infty_{\tau}|\lambda^j + |\tilde{x}_{\tau+N}^*-x^\infty_{\tau+N}|\lambda^{N-j}\right)
	\end{equation*}
	for all $j\in[0,N]$, $\tau\in\mathbb{N}$, $\bar{x}\in\mathcal{X}$, and all possible $d_{-\infty:\infty}$.
\end{ass}

Assumption~\ref{ass:turnpike_inf_prior} essentially states that the infinite-horizon solution $z^\infty$ serves as turnpike for MHE problems with prior weighting, cf.~Section~\ref{sec:turnpike_sensitivity}. Note that such behavior could be observed in all numerical examples in Section~\ref{sec:ex:MHE}.

\begin{rem}[Exponential turnpike]
	In contrast to Definition~\ref{def:turnpike}, we impose an exponential turnpike property in Assumption~\ref{ass:turnpike_inf_prior}, which is crucially required to derive uniform (exponential) convergence of the turnpike prior $\bar{x}_{t-N_t}$ to the turnpike $x^\infty$ in the following proposition. Note that this is conceptually similar to recent stability results for nonlinear MHE, which also require exponential detectability (rather than asymptotic detectability), cf., e.g., \cite{Rawlings2017,Allan2021a,Schiller2023c,Hu2023}.
\end{rem}

\begin{prop}\label{prop:prior}
	Let Assumption~\ref{ass:turnpike_inf_prior} hold and the sets $\mathcal{X},\mathcal{W},\mathcal{V}$ be compact.
	Suppose there exists $p_1,p_2>0$ and $a\geq1$ such that
	\begin{equation}\label{eq:bound_prior}
		p_1|x-\bar{x}|^a \leq p_t(x,\bar{x}) \leq p_2|x-\bar{x}|^a
	\end{equation}
	for all $x,\bar{x}\in\mathcal{X}$ uniformly for all $t\in\mathbb{N}$. Furthermore, assume that there exists $\alpha_1,\alpha_2,\alpha_3\in\mathcal{K}$ such that
	\begin{align}
		&l(\hat{x},\hat{w};(u,y)) \leq \alpha_1(|\hat{w}|) + \alpha_2(|y-h(\hat{x},u)|),\label{eq:bound_l}\\
		&g(\hat{x};(u,y)) \leq  \alpha_3(|y-h(\hat{x},u)|)\label{eq:bound_g}
	\end{align}
	for all $\hat{x}\in\mathcal{X}$, $\hat{w}\in\mathcal{W}$, $(u,y)\in\mathcal{U}\times\mathbb{R}^p$ satisfying $y-h(\hat{x},u)\in\mathcal{V}$.	
	Then, there exists $\sigma\in\mathcal{L}$ such that for any $\rho\in(0,1)$, there exists $\bar{N}\in\mathbb{N}_\mathrm{e}$ such that the turnpike prior~\eqref{eq:prior_turnpike} satisfies
	\begin{align}\label{eq:prop:turnpike}
		|\bar{x}_{t-N/2} - {x}^\infty_{t-N/2}| \leq \rho |\bar{x}_{t-N}-x^\infty_{t-N}| + \sigma(N)
	\end{align}
	for all $N\in\mathbb{N}_\mathrm{e}$, $N\geq\bar{N}$ and all $t\in[N,\infty)$.
\end{prop}

\begin{proof}
	Consider some data sequence $d_{-\infty:\infty}$ associated with the system~\eqref{eq:sys} and the respective infinite-horizon solution~$z^\infty$.
	Let $t\in[N,\infty)$ and $\bar{x}_{t-N}\in\mathcal{X}$ be arbitrary and denote the solution of $P^\mathrm{p}_N(d_{t-N:t},\bar{x}_{t-N},t)$ by $\tilde{z}^*_{t-N:t}$.
	From Assumption~\ref{ass:turnpike_inf_prior}, we obtain
	\begin{equation}
		|\tilde{z}^*_{t-N+j}-z^\infty_{t-N+j}| \leq K|\tilde{x}^*_{t-N}-x^\infty_{t-N}|\lambda^{j} + Kc_1\lambda^{N-j} \label{eq:proof_prior_0}
	\end{equation}
	for all $j\in[0,N]$, where $c_1>0$ satisfies $|x_1-x_2|\leq c_1$ for all $x_1,x_2\in\mathcal{X}$ (compactness of $\mathcal{X}$ ensures existence of such $c_1$).
	By the triangle inequality, it follows that
	\begin{equation}
		|\tilde{x}_{t-N}^*-x^\infty_{t-N}| \leq |\tilde{x}_{t-N}^*-\bar{x}_{t-N}| + |\bar{x}_{t-N}-x^\infty_{t-N}|. \label{eq:proof_prior_1}
	\end{equation}
	Using \eqref{eq:bound_prior}, the cost function~\eqref{eq:J_p}, and optimality of $\tilde{z}_{t-N:t}^*$, it holds that
	\begin{align}
		&p_1|\tilde{x}_{t-N}^*-\bar{x}_{t-N}|^a \leq J^\mathrm{p}_N(\tilde{x}_{t-N:t}^*,\tilde{w}_{t-N:t-1}^*;d_{t-N:t},t)\nonumber \\
		&\leq J^\mathrm{p}_N({x}^\infty_{t-N:t},{w}^\infty_{t-N:t-1};d_{t-N:t},t)\nonumber\\
		&\leq p_2|{x}^\infty_{t-N}{-\,}\bar{x}_{t-N}|^a{+\,} J_{N}({x}^\infty_{t-N:t},{w}^\infty_{t-N:t-1};d_{t-N:t}). \label{eq:proof_prior_2}
	\end{align}
	The bounds in \eqref{eq:bound_l} and \eqref{eq:bound_g}, compactness of $\mathcal{W}$ and $\mathcal{V}$, and similar arguments as in the proof of Lemma~\ref{lem:VT} imply the existence of uniform constants $A,B>0$ such that $J_{N}({x}^\infty_{t-N:t},{w}^\infty_{t-N:t-1};d_{t-N:t}) \leq AN+B$. In combination with \eqref{eq:proof_prior_1} and \eqref{eq:proof_prior_2} and using the fact that the function $s\mapsto s^{1/a}$ is subadditive for $s\geq0$ and $a\geq 1$, this leads to
	\begin{equation*}
		|\tilde{x}_{t-N}^*-x^\infty_{t-N}|\leq c_2|{x}^\infty_{t-N}-\bar{x}_{t-N}| + \left(\frac{AN+B}{p_1}\right)^{1/a},
	\end{equation*}
	where $c_2:=1+({p_2}/{p_1})^{1/a}$. Evaluating \eqref{eq:proof_prior_0} at $j=N/2$, we can infer that
	\begin{align}
		|\tilde{z}_{t-N/2}^*-z^\infty_{t-N/2}|&\leq Kc_2|{x}^\infty_{t-N}-\bar{x}_{t-N}|\lambda^{N/2}\label{eq:proof_prior_3}\\
		& \quad + K\left(\frac{AN+B}{p_1}\right)^{1/a}\lambda^{N/2}+Kc_1\lambda^{N/2}. \nonumber
	\end{align}
	For any $\rho\in(0,1)$, there exists $\bar{N}\in\mathbb{N}_\mathrm{e}$ sufficiently large such that $Kc_2\lambda^{N/2}\leq\rho$ for all $N\geq\bar{N}$.
	Define $\sigma_1(r) := K((Ar+B)/p_1)^{1/a}\lambda^{r/2}$, $r\geq0$.
	Note that $\sigma_1(r)$ converges to zero for $r\rightarrow\infty$, since the exponential term dominates for large enough $r$.
	Hence, there exists $\bar{\sigma}_1\in\mathcal{L}$ satisfying $\sigma_1(r)\leq \bar{\sigma}_1(r)$ for all $r\geq0$. Thus, from~\eqref{eq:proof_prior_3}, we can infer that the updated turnpike prior $\bar{x}_{t-N/2} = \tilde{x}_{t-N/2}^*$ satisfies
	\begin{align*}
		|\bar{x}_{t-N/2}-x^\infty_{t-N/2}| &\leq |\tilde{z}^*_{t-N/2}-z^\infty_{t-N/2}| \\
		&\leq \rho|\bar{x}_{t-N}-x^\infty_{t-N}|  + \bar{\sigma}_1(N) + Kc_1\lambda^{N/2}
	\end{align*}
	for all $t\geq N$ and all $N\in\mathbb{N}_\mathrm{e}$ satisfying $N\geq\bar{N}$. Defining $\sigma(r):=\bar{\sigma}_1(r)+Kc_1\lambda^{r/2}$, $r\geq0$ and noting that $\sigma\in\mathcal{L}$ establishes the statement and hence finishes this proof.
\end{proof}

The conditions in \eqref{eq:bound_prior}--\eqref{eq:bound_g} on the prior weighting, stage cost, and terminal cost are standard (cf., e.g., \cite[Ass.~4.22]{Rawlings2017}) and obviously satisfied for the practically relevant case of quadratic penalties as in~\eqref{eq:stage_cost}, \eqref{eq:term_cost}, and \eqref{eq:prior_cost}.
Provided that the horizon length $N$ is chosen sufficiently large, Proposition~\ref{prop:prior} implies (by a recursive application of the bound in \eqref{eq:prop:turnpike}) that the prior $\bar{x}_{t-N_t}$~\eqref{eq:prior_turnpike} forms a sequence that exponentially converges into a neighborhood of the turnpike (the infinite-horizon benchmark solution $x^\infty$), which is also evident in the simulation example in Section~\ref{sec:ex:MHE_CSTR}. Here, we want to emphasize that the size of this neighborhood depends on the horizon length~$N$ and can in fact be made arbitrarily small by choosing larger values of $N$ (due to the fact that $\sigma\in\mathcal{L}$).

Overall, a properly selected prior weighting $p_t$ according to Proposition~\ref{prop:prior} ensures that the initial state $\tilde{x}_{t-N_t|t}^*$ of the solution of the finite-horizon problem $P_{N_t}^\mathrm{p}(\cdot,t)$ is close to the turnpike, which hence effectively reduces the approaching arc and allows using short horizons. However, the natural occurrence of the leaving arc can still cause the resulting estimate $\hat{x}_t^{\mathrm{mhe,p}} = \tilde{x}^*_{t|t}$ to be again relatively far away from the turnpike, which could again be reduced by using an artificial delay in the MHE scheme as suggested in Section~\ref{sec:MHE:delay}.

\begin{rem}[Performance of $\delta$MHE with prior weighting]\label{rem:perf}
	For some fixed delay $\delta\in[0,N/2]$, we can define $\delta$MHE (with prior weighting) as
	\begin{equation}\label{eq:MHE_delay_x_prior}
		\hat{x}_{t-\delta}^{\mathrm{\delta mhe,p}} =  \tilde{x}_{t-\delta|t}^*, \ \ t\in[ \delta,\infty).
	\end{equation}
	Under Assumption~\ref{ass:turnpike_inf_prior}, it is straightforward to show that $|\hat{x}^{\delta \mathrm{mhe,p}}_{j}-x^\infty_{j}|\leq \sigma(\delta)$ for all $j\in[\delta,t-\delta]$, $t\geq\delta$, with $\sigma\in\mathcal{L}$ from Proposition~\ref{prop:dMHE} (which can be easily modified to this case).
	As a result, the performance estimates from Theorem~\ref{thm:perf_MHE}, Corollary~\ref{cor:regret}, and Corollary~\ref{cor:performance} directly carry over to $\delta$MHE with prior weighting.
	Thus, $\delta$MHE with prior weighting and a suitably selected delay $\delta$ can recover the accuracy and performance of the infinite-horizon estimator, with shorter horizons compared to $\delta$MHE without prior weighting.
	Here, we want to emphasize that this conclusion holds under Assumption~\ref{ass:turnpike_inf_prior}, i.e., for any choice of the prior estimate $\bar{x}_{t-N_t}$ from~\eqref{eq:prior_filtering}--\eqref{eq:prior_turnpike} (in contrast to Proposition~\ref{prop:prior}, which is an exclusive feature of the proposed turnpike prior~\eqref{eq:prior_turnpike}).
\end{rem}

\section{Batch state estimation}\label{sec:offline}

In this section, we consider the case of batch (i.e., offline) state estimation.
Here, one is interested in matching a given data batch $d_{0:T}$ for some $T\in\mathbb{N}$ to the system equations~\eqref{eq:sys} to recover the true unknown state sequence $x_{0:T}$.
This is particularly relevant for large-scale nonlinear data assimilation problems in geophysics and environmental sciences, such as numerical weather prediction, oceanography, and hydrology; see, e.g., \cite{Asch2016,Carrassi2018}.
To this end, a natural approach is trying to compute $\hat{x}^*_{0:T}$ by solving the optimal state estimation problem in~\eqref{eq:MHE} with $N=T$.
However, if the data set (i.e., $T$) or the underlying model is large, or if the computations are limited in terms of time or resources, solving the full problem $P_T(d_{0:T})$ for the optimal solution becomes difficult, if not impossible.

Instead, we approximate the optimal state sequence $\hat{x}_{0:T}^*$ by solving a sequence of truncated problems $P_N$ of length $N\in\mathbb{N}_\mathrm{e}$ with $N\ll T$ on overlapping data windows.
Based on the results of Section~\ref{sec:online}, only the central segment of each solution corresponding to the interval $N/2\pm\Delta$ for some $\Delta\in[0,N/2]$ is retained (where the turnpike property ensures high accuracy), while boundary regions (approaching and leaving arcs) are discarded.
The retained segments are then concatenated to cover the whole interval $[0,T]$.

\begin{figure}[t]
	\centering
	\includegraphics[width=0.97\columnwidth]{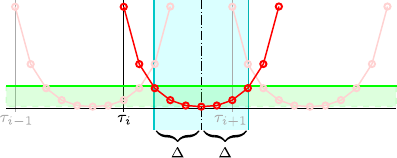}
	\caption{Construction of the approximate estimator. Depicted is the difference $|\zeta_N^x(j,d_{\tau_i:\tau_i+N})-\hat{x}^*_{\tau_i+j}|$ over $j\in[0,N]$ for a fixed $\tau_i$ (red), along with the central interval $N/2\pm \Delta$ (cyan); the green-colored region includes the solution segments used to form the approximate estimator in~\eqref{eq:candidate_x}, with contributions from neighboring windows at $\tau_{i-1}$ and $\tau_{i+1}$ (light red).}
	\label{fig:AE}
\end{figure}

For simplicity, assume hat there exists some $k\in\mathbb{N}$ such that $T$, $N$, and $\Delta$ satisfy $T=N+(k+1)(2\Delta+1)$.
Then, the proposed approximation requires solving $K_{\Delta} := k+2=({T-N})/({2\Delta+1})+1$ truncated problems $P_N(d_{\tau_i:\tau_i+N})$ with $\tau_i=i(2\Delta+1)$, $i\in[0,K_\Delta-1]$.
This leads to the definition of the \emph{approximate estimator} $\hat{x}^\mathrm{ae}_t$ for $t\in[0,T]$ as
\begin{equation}
	\hat{x}^\mathrm{ae}_{\tau_i+j} {\,=}
	\begin{cases}
		\zeta_N^x(j,d_{0:N}), j \in[0,N/2+\Delta],
		& \hspace{-1ex} i=1\\
		\zeta_N^x(j,d_{\tau_i:\tau_i+N}), j{\,\in\,} N/2\pm \Delta,
		& \hspace{-1ex} i{\,\in\,}[1,K_\Delta{-\,}2]\\
		\zeta_N^x(j,d_{T-N:T}), j{\,\in\,}[N/2-\Delta,N],
		& \hspace{-1ex} i=K_\Delta-1,
	\end{cases}\label{eq:candidate_x}
\end{equation}
where we use the mapping $\zeta_N^x$ as defined below~\eqref{eq:z_N}. The construction is also illustrated in Figure~\ref{fig:AE}.

For $\Delta = 0$, the approximate estimator in~\eqref{eq:candidate_x} corresponds to the $\delta$MHE scheme in~\eqref{eq:MHE_delay_x} with delay $\delta=N/2$.
Hence, the accuracy and performance results from Section~\ref{sec:online} directly apply with respect to the benchmark~$x^\infty$ (where $d_{-\infty:\infty}$ denotes a suitable extension of $d_{0:T}$ satisfying~\eqref{eq:sys}--\eqref{eq:con_Z} for all $t\in\mathbb{Z}$).
Therefore, under Assumption~\ref{ass:turnpike_inf}, the estimates $\hat{x}^\mathrm{ae}_j$ are close to the turnpike $x^\infty_j$ for all $j\in[N/2,T-N/2]$.
The same property holds for the desired optimal solution $\hat{x}^*_{0:T}$ of $P_T(d_{0:T})$, implying that both sequences are close to $x^\infty$ (and hence to each other) on this interval, compare~Figure~\ref{fig:turnpike_sol}. Consequently, $\hat{x}^\mathrm{ae}_{0:T}$ is approximately optimal on $[N/2, T{-}N/2]$ with respect to $\hat{x}^*_{0:T}$.
For $\Delta{\,\in\,}[1,N/2]$, analogous performance guarantees hold, albeit with bounds that degrade as $\Delta$ increases.

\begin{rem}[Computational{\,}effort{\;}vs{.\,}accuracy]
	\hspace{-0.5ex}The parameter $\Delta$ determines the number of retained elements and thus the trade-off between computational effort and accuracy; increasing $\Delta$ decreases $K_\Delta$ while increasing the approximation error, and \textit{vice versa} (see Figure~\ref{fig:AE}).
	Under an exponential turnpike property, it suffices to choose $\Delta$ close to $N/2$, so that only a few elements at the beginning and end of each truncated solution---corresponding to the approaching and leaving arcs---are discarded. This can lead to a significant reduction in computational effort while incurring only negligible performance degradation; see Section~\ref{sec:ex:OE_LTI}.
\end{rem}

\begin{rem}[Distributed optimization]
	The approximation proposed in~\eqref{eq:candidate_x} decomposes the full problem into a sequence of smaller, independent subproblems, thereby accelerating computation in two ways: (i) each subproblem is significantly smaller than the full problem, and (ii) the subproblems are fully decoupled and can be solved in parallel.
	At the same time, near-optimal performance is preserved, as the turnpike property ensures that the central segments of the truncated solutions closely approximate the desired optimal solution.
	Consequently, \eqref{eq:candidate_x} provides a distributed optimization method for the full problem, where the error can be made arbitrarily small by tuning $N$ and $\Delta$, as illustrated in Section~\ref{sec:ex:OE_LTI}.
	This is particularly relevant for large-scale data assimilation problems~\cite{Asch2016,Carrassi2018}, and more broadly in the context of robust optimization for data-driven decision making; cf.~\cite{MohajerinEsfahani2018,MohajerinEsfahani2017a}.
\end{rem}

\section{High performance implies high accuracy}\label{sec:accuracy}
In this section, we relate the performance of a state estimator, measured by~\eqref{eq:perf_criterion}, to its accuracy in terms of the estimation error. To this end, we impose a detectability condition to ensure that the measurement data contains sufficient information about the unknown state trajectory, where we use the notion of incremental input/output-to-state stability (i-IOSS).

\begin{ass}[Exponential i-IOSS]\label{ass:IOSS}
	The system~\eqref{eq:sys} is exponentially i-IOSS, i.e., there exists a continuous function $U:\mathcal{X}\times\mathcal{X}\rightarrow\mathbb{R}_{\geq 0}$ together with matrices $P_1,P_2,Q,R\succ0$ and a constant $\eta\in[0,1)$ such that
	\begin{subequations}
		\label{eq:IOSS}
		\begin{align}\label{eq:IOSS_bounds}
			&|x_1-x_2|_{P_1}^2 \leq U(x_1,x_2) \leq |x_1-{x}_2|_{P_2}^2, \\[1ex]
			&U(f(x_1,u,w_1),f(x_2,u,w_2))\label{eq:IOSS_supply}\\
			& \leq \eta U(x_1,x_2) +|w_1-w_2|_Q^2 + |h(x_1,u)-h(x_2,u)|_R^2 \nonumber
		\end{align}
	\end{subequations}
	for all $(x_1,u,w_1),(x_2,u,w_2)\in\mathcal{X}\times\mathcal{U}\times\mathcal{W}$.
\end{ass}

Assumption~\eqref{ass:IOSS} is a Lyapunov characterization of exponential i-IOSS, which became a standard detectability condition in the context of MHE in recent years, see, e.g., \cite{Allan2021a,Hu2023,Knuefer2023,Rawlings2017,Schiller2023c,Alessandri2025}. This property implies that the difference between any two state trajectories is bounded by the differences of their initial states, their disturbance inputs, and their outputs. In other words, if two trajectories of an i-IOSS system have the same inputs and outputs, then their states must converge to each other.
Note that Assumption~\ref{ass:IOSS} is not restrictive in the state estimation context; in fact, by adapting the results from \cite{Allan2021,Knuefer2023}, it is in fact necessary and sufficient for the existence of robustly (exponentially) stable state estimators. Moreover, Assumption~\ref{ass:IOSS} can be verified using LMIs, cf., e.g, \cite{Schiller2023c,Arezki2023}.

Given a sequence $\hat{x}_{0:t}$, $t\in\mathbb{N}$ produced by some state estimator, the following result establishes a bound with respect to the true state sequence ${x}_{0:t}$ in terms of its performance.

\begin{prop}\label{prop:IOSS}
	Suppose Assumption~\ref{ass:IOSS} holds.
	Consider the performance measure~\eqref{eq:perf_criterion} for some $t_1,t_2\in\mathbb{N}$ with the quadratic stage cost $l$ from~\eqref{eq:stage_cost} for some $Q,R\succ0$.
	Then, there exist $C_1,C_2,C_3>0$ such that
	\begin{align}
		&|{x}_\tau-\hat{x}_\tau|^2 \leq  C_1\eta^{\tau-t_1} |{x}_{t_1}-\hat{x}_{t_1}|^2 \label{eq:prop_accuracy} \\
		&\hspace{1.9cm} + C_2\max_{j\in[t_1,\tau-1]}\{|w_j|^2,|v_j|^2\} + C_3J_{[t_1,t_2]}(\hat{z}_{t_1:t_2}) \nonumber
	\end{align}
	for all $\tau\in[t_1,t_2]$, all initial conditions $x_0,\hat{x}_0\in\mathcal{X}$, and all input and disturbance sequences $u\in\mathcal{U}^\infty$, $w,\hat{w}\in\mathcal{W}^\infty$, and $v\in\mathcal{V}^\infty$, where $\eta\in(0,1)$ is from Assumption~\ref{ass:IOSS} and $x_{j+1}=f(x_{j},u_{j},w_{j})$, $y_{j} = h(x_{j},u_{j}) + v_{j}$, $\hat{x}_{j+1}=f(\hat{x}_{j},u_{j},\hat{w}_{j})$, and $\hat{z}_{j} = (\hat{x}_{j},\hat{w}_{j})$ for all $j\in\mathbb{N}$.
\end{prop}

\begin{proof}
	Without loss of generality, we can assume that the matrices $Q$ and $R$ from the cost function~\eqref{eq:stage_cost} are the same\footnote{
		For any weighting matrices $Q,R\succ0$, the i-IOSS Lyapunov function $U$ can be suitably scaled such that \eqref{eq:IOSS_supply} holds with that choice of $Q$ and $R$.
	} matrices as in~\eqref{eq:IOSS_supply}.
	For any $\tau\in[t_1,t_2]$, the sequences $x_{0:\tau}$ and $\hat{x}_{0:\tau}$ form trajectories of the i-IOSS system~\eqref{eq:sys}.
	Using the dissipation inequality~\eqref{eq:IOSS_supply} from Assumption~\ref{ass:IOSS}, the fact that $|a-b|_H^2\leq 2|a|_H^2 + 2|b|_H^2$ for any real vectors $a,b$ and matrix $H\succ0$ by Cauchy-Schwarz and Young's inequality, and the performance criterion in~\eqref{eq:perf_criterion}, we can infer that
	\begin{align*}
		&U({x}_\tau,\hat{x}_\tau)
		\leq \eta^{\tau-t_1} U({x}_{t_1},\hat{x}_{t_1}) \\
		&\hspace{1.2cm} + 2\sum_{j=t_1}^{\tau-1}\eta^{\tau-j-1}(|{w}_{j}|_Q^2 + |v_{j}|_R^2) + 2J_{[t_1,t_2]}(\hat{z}_{t_1:t_2}).
	\end{align*}
	The fact that $a+b\leq \max\{2a,2b\}$ for all $a,b\geq0$ together with the convergence property of the geometric series leads to
	\begin{equation*}
		\sum_{j=t_1}^{\tau-1}\eta^{\tau-j-1}(|{w}_{j}|_Q^2 + |v_{j}|_R^2)
		\leq \frac{2C}{1{\,-\,}\eta}\max_{j\in[t_1,\tau-1]}\{|w_j|^2,|v_j|^2\}
	\end{equation*}
	with $C:=\max\{\lambda_{\max}(Q),\lambda_{\max}(R)\}$.
	Using \eqref{eq:IOSS_bounds}, we obtain~\eqref{eq:prop_accuracy} with $C_1 = \frac{\lambda_{\max}(P_2)}{\lambda_{\min}(P_1)}$, $C_2 = 4C/(\lambda_{\min}(P_1)(1-\eta))$, and $C_3= 2/\lambda_{\min}(P_1)$, which finishes this proof.	
\end{proof}

Proposition~\ref{prop:IOSS} draws a link between the performance of a state estimator (measured by the criterion in~\eqref{eq:perf_criterion}) and the corresponding estimation error.
In particular, for large $\tau$, the error is upper bounded by the performance $J_{[t_1,t_2]}$ and the maximum disturbance $w_j,v_j$, $j\in[t_1,\tau-1]$ that affected the past system behavior and associated measurement data.
Hence, if the disturbances are small, we directly have that good performance (small values of $J_{[t_1,t_2]}$) implies accurate estimation results (small errors $|x_\tau-\hat{x}_\tau|$).
Consequently, it is indeed advisable to design estimators that achieve good performance in the sense of the criterion $J_{[t_1,t_2]}$.

We want to emphasize that Proposition~\ref{prop:IOSS} and its implications apply for \emph{any} state estimator/observer design.
However, as the performance criterion appearing in the estimation error bound in \eqref{eq:prop_accuracy} constitutes a part of the cost function used in the infinite-horizon problem $P_\infty$, the corresponding solution $z^\infty$ provides a comparatively small bound for the estimation accuracy for all $\tau\in[t_1,t_2]$ for \emph{any} $t_1,t_2\in\mathbb{N}$.
Under the turnpike condition from Assumption~\ref{ass:turnpike_inf}, we know that $\delta$MHE along with a suitably selected delay $\delta$ achieves nearly the same performance as the benchmark $z^\infty$ on any interval $[t_1,t_2]\subseteq[\delta,t-\delta]$ for all $t\geq\delta$, and is hence expected to be similarly\footnote{
	This also applies to $\delta$MHE with prior weighting under Assumption~\ref{ass:turnpike_inf_prior}.
} accurate.
This highly useful feature of $\delta$MHE can also be seen in all our simulation examples in Section~\ref{sec:example}.

\section{Numerical examples}\label{sec:example}

We illustrate the theoretical results from the previous sections by means of several simulation\footnote{
	{The code is publicly available in the GitHub repository at \url{https://github.com/jdschille/mhe-turnpikes} \cite{Schiller2025_code}.}
	} examples.
In Section~\ref{sec:ex:OE_LTI}, we focus on batch state estimation, where we consider linear time-invariant (LTI) systems with up to 200 states.
This allows us to compute the globally optimal solution to the full state estimation problem as well as to apply standard optimal filtering and smoothing techniques for comparison reasons.
In line with our theoretical findings from Section~\ref{sec:offline}, the proposed approximate estimator (AE) approximates the desired optimal solution with negligible error, while reducing the required computation time by more than 95\%. In Section~\ref{sec:ex:MHE}, we focus on online state estimation for nonlinear systems (in particular, MHE with prior weighting) and consider two realistic examples from the literature: a continuous stirred-tank reactor and a highly nonlinear 12-state quadrotor model.
In both examples, we observe the turnpike behavior.
Our main finding is that already a small delay in the MHE scheme (one to three steps) reduces the overall estimation error with respect to the true unknown system state by 20\,--25\%.

\subsection{Batch state estimation}\label{sec:ex:OE_LTI}

{We apply the proposed AE to large batch estimation problems, where the direct computation of the optimal solution is either time-consuming or simply impossible.
Specifically, we consider random observable, stable LTI systems of the form
\begin{equation*}
	x^+ = Ax + Bu + w, \quad 
	y = Cx + v,
\end{equation*}
where we select $A$, $B$, $C$ using the \texttt{drss} command of MATLAB with $m=30$ inputs, $n\in\{50,100,200\}$ states, and $p\in\{20,40,80\}$ outputs.}
{We consider a given data batch $d_{0:T}$ with $T=4803$, where the system was subject to $x_0=0$, a known sinusoidal input sequence $u_{0:T}$, and unknown process disturbance $w_{0:T-1}$ and measurement noise~$v_{0:T}$ (during simulations, $w_j$ is sampled from a Laplace distribution superimposed with a deterministic sinusoidal function that may represent unmodeled nonlinear dynamics, and $v_j$ is sampled from a normal distribution superimposed with random impulsive noise affecting approximately 2\% of the samples).}

\begin{figure}[t]
	\centering
	\includegraphics[scale=1]{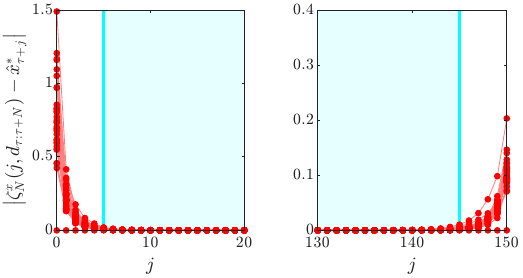}
	\caption{Differences between the optimal solution $\hat{x}^*_{\tau_i+j}$ solving the full problem $P_T(d_{0:T})$ and the truncated solutions required to construct the proposed AE, depicted over $j\in[0,N]$ for $\tau_i=i(2\Delta+1)$, $i=0,...,K_\Delta-1$. The cyan-colored region includes all solution elements corresponding to $j\in N/2\pm\Delta$ that are used to form the approximate estimator.}
	\label{fig:LTI_TP}
\end{figure}

For reconstructing the unknown state sequence $x_{0:T}$, we choose the quadratic cost function~\eqref{eq:MHE_cost}--\eqref{eq:term_cost} and select $Q=I_n$ and $R=G=I_p$.
In the following, we compare the optimal solution to the full problem $P_T(d_{0:T})$ with the proposed AE from~\eqref{eq:candidate_x}, where we select $k=32$, $N=150$, and $\Delta=70$.
We hence only discard the first and last five elements of each truncated solution involved in computing the AE.
As Figure~\ref{fig:LTI_TP} shows, this is already sufficient to effectively remove the influence of the approaching and leaving arcs, mainly due to the fact that the optimal state estimation problem exhibits \emph{exponential} turnpike behavior, compare Section~\ref{sec:turnpike_sensitivity}.

The full and truncated optimal estimation problems can be cast as unconstrained quadratic programs (QPs), which we solve using the \texttt{quadprog} implementation of MATLAB. The truncated QPs for the AE are additionally solved in parallel using the \textit{Parallel Computing Toolbox}.
{For comparison reasons, we additionally consider the Kalman filter (KF) using the covariance matrices $Q^{-1}$ and $R^{-1}$, where the initial estimate $\hat{x}^{\mathrm{kf}}_0$ is drawn from an isotropic normal distribution with zero mean and covariance matrix $\hat{P}_0^{\mathrm{kf}} = 10^3\cdot I_n$.
Moreover, we also apply the fixed-interval Rauch-Tung-Striebel (RTS) smoother over the entire data batch, which provides the best possible (i.e., optimal) estimate in the context of linear Gaussian models, see \cite[Ch.~5]{Crassidis2011} and \cite[Ch.~8]{Saerkkae2013} for further details and a description of the corresponding algorithm.}

Table~\ref{tab:J_LTI} shows the estimation results in terms of the performance index $J_{[0,T]}$, accuracy (SSE), and overall computation time~$t_\mathrm{comp}$ for different system dimensions. For $n\in\{50,100\}$, we observe that the performance $J_{[0,T]}$ and the SSE of the full solution and the AE are nearly identical. However, the AE can be computed more than 20 times faster than the full solution (saving more than 95\% of the computation time) due to a (much) smaller problem size and the fact that the truncated problems are solved in parallel.
For $n=120$, it was already impossible to numerically solve the full problem $P_T(d_{0:T})$ due to the problem size (in contrast to the AE).
For all system realizations, the KF is much faster as expected, because the QPs are replaced by simple matrix computations---however, it performs much worse than the AE (both in terms of $J_{[0,T]}$ and SSE), which is mainly due to the fact that only past data is used to compute the corresponding estimates.
The RTS combines the KF forward recursion with a backward pass, so that each estimate exploits the full data batch, at the cost of increased computations. While this improves the estimation results, both performance and SSE remain inferior to the AE. This is mainly due to the necessary use of the initial estimate in the recursions, which is avoided in the cost function~\eqref{eq:MHE_cost}.

\begin{table}[t]
	\centering
	\caption{Estimation results of the proposed approximate estimator (AE) compared to the solution of the optimal state estimation problem $P_T(d_{0:T})$ (optimal), the Kalman filter (KF), and the fixed-interval RTS smoother (RTS). The parentheses contain relative changes with respect to the optimal solution.}
	\label{tab:J_LTI}
	\setlength\tabcolsep{3.8pt}
	\begin{tabular}{ccccc}
		\toprule
		($n$,$p$) & Method  & $J_{[0,T]}$ & $\mathrm{SSE}$ & $t_\mathrm{comp}\, [\mathrm{s}]$ \\\midrule
		\multirow{4}{*}{\rotatebox{90}{(50,20)}}
		& optimal & 62.646 &   167.41 &   20.472\\ 
		& AE&   62.651 (+0.007\%) & 167.42 (+0.003\%) &    0.927 (-95.5\%)\\ 
		& RTS&  63.271 (+1\%) \hspace{3ex}  &   441.61 (+163\%) \hspace{0.75ex} & 0.868 (-95.8\%)\\
		& KF& 4.9·10\textsuperscript{5} (+777\%) \hspace{1ex} & 3.3·10\textsuperscript{4} (+198\%) \hspace{1ex} &    0.374 (-98.2\%)  \\ 
		\midrule
		\multirow{4}{*}{\rotatebox{90}{(100,40)}}
		& optimal& 83.667 &   205.33 &   80.845\\
		& AE&   83.675 (+0.009\%) &   205.35 (+0.009\%) &    3.809 (-95.3\%)\\ 
		& RTS&   85.754 (+2\%) \hspace{1.7ex} \phantom{,}&   634.03 (+209\%) \hspace{-0.6ex} \phantom{,}&    2.364 (-97.1\%)\\ 
		& KF& 4·10\textsuperscript{6}& 6·10\textsuperscript{4} &    1.158  (-98.6\%) \\ 
		\midrule
		\multirow{3}{*}{\rotatebox{90}{(200,80)}}
		& AE&  123.362 &   276.90 &   20.133\\ 
		& RTS&  125.225 &   749.22 &   12.729\\ 
		& KF& 5.2·10\textsuperscript{8} & 9.9·10\textsuperscript{4} &    4.896\\ 
		\bottomrule
	\end{tabular}
\end{table} 

Overall, we find that the AE approximates the desired optimal solution with negligible error, which is particularly important in practice when the full problem cannot be directly solved due to its size or complexity and iterative solutions such as the KF and related smoothing algorithms are not sufficiently accurate.

\subsection{Online estimation}\label{sec:ex:MHE}

{We now focus on the case of online state estimation.}
{The following simulations were performed on a standard laptop in MATLAB using CasADi \cite{Andersson2018} and the solver IPOPT \cite{Waechter2005}.}

\subsubsection{Continuous stirred-tank reactor}\label{sec:ex:MHE_CSTR}
We consider the continuous stirred-tank reactor (CSTR) from \cite[Example 1.11]{Rawlings2017}, where an irreversible, first order reaction $A\rightarrow B$ occurs in the liquid phase and the reactor temperature is regulated with external cooling, see also~\cite{Pannocchia2003} for more details. The continuous-time nonlinear state space model is given by
\begin{align*}
	&\frac{dc}{dt} = \frac{F_0(c_0-c)}{\pi r^2h} - k_0\exp\left(\frac{-E}{RT}\right)c,\\
	&\frac{dT}{dt}{\,=\,}\frac{F_0(T_0{\,-}T)}{\pi r^2h} {\,+\,} \frac{-\Delta H}{\rho C_p}k_0\exp\left(\frac{-E}{RT}\right)c {\,+\,} \frac{2U}{r\rho C_p}(T_c{-\,}T),\\
	&\frac{dh}{dt} = \frac{F_0-F}{\pi r^2},
\end{align*}
where the states are $c$ (the molar concentration of species A), $T$ (the reactor temperature), and $h$ (the level of the tank), and the control inputs are $T_c$ (the coolant liquid temperature) and $F$ (the outlet flowrate).
The model parameters are taken from \cite[Example 1.11]{Rawlings2017}.
The open-loop stable steady state solution is $x^\mathrm{ss}=[0.878,323.5,0.659]^\top$ associated with the input $u^\mathrm{ss} = [300,0.1]^\top$.
By using Runge-Kutta discretization (RK4) with sampling time $t_\Delta=0.25$, we obtain the discrete-time model in~\eqref{eq:sys} with $f$ as in \eqref{eq:dynamics} and $h(x,u) = x_2$, where we define $x = [c,T,h]^\top$ and $u = [T_c,F]^\top$, consider additional process disturbance $w\in\mathbb{R}^3$, and assume that only the temperature of the reactor $T$ can be measured, subject to the measurement noise $v\in\mathbb{R}$.
The real system is initialized with $x_0=[0.8, 295, 0.7]^\top$. In the following, we consider the simulation time $T_\mathrm{s} = 200$ and apply an open-loop control sequence $u_{0:T_\mathrm{s}}$ with $u_2 = u^\mathrm{ss}_2$ and $u_1$ as a trapezoidal input sequence with plateaus at $u^{\mathrm{ss}}_1$ and $u^{\mathrm{ss}}_1-25$.
During the simulations, we sample the disturbances $w$ and $v$ from uniform distributions over $\{w: |w_{1}|\leq5\cdot10^{-3}, |w_2|\leq 1, |w_{3}|\leq5\cdot10^{-3}\}$ and $\{v: |v|\leq 3\}$.
To estimate the true unknown state $x_t$ from the measurement data, we design different MHE schemes that rely on the cost function~\eqref{eq:J_p} with the horizon length $N=10$ and quadratic costs~\eqref{eq:stage_cost} and \eqref{eq:term_cost}, where we select $Q = \mathrm{diag}(10^3,1,10^5)$, $R=G=1$.
Moreover, we use a quadratic prior weighting \eqref{eq:prior_cost} with time-varying matrix $W_t$, $t\in\mathbb{N}$, initialized with $W_0=10^{-2}I_3$ and updated using the well-known covariance formulas of the extended Kalman Filter (EKF).
In the optimal estimation problems, we also consider the state constraints $\mathcal{X}=[0.5,1.2]\times[200,400]\times[0.5,1.2]$, but we consider the disturbance sets to be unknown and use $\mathcal{W}=\mathbb{R}^3$ and $\mathcal{V}\in\mathbb{R}$.
In the following, we compare MHE with filtering prior~\eqref{eq:prior_filtering}, smoothing prior~\eqref{eq:prior_smoothing}, and turnpike prior~\eqref{eq:prior_turnpike} and additionally consider the infinite-horizon estimator (IHE), which we approximate by solving the clairvoyant FIE problem using all available data $d_{0:T_\mathrm{s}}$.
\begin{figure}[t]
	\centering
	\includegraphics[]{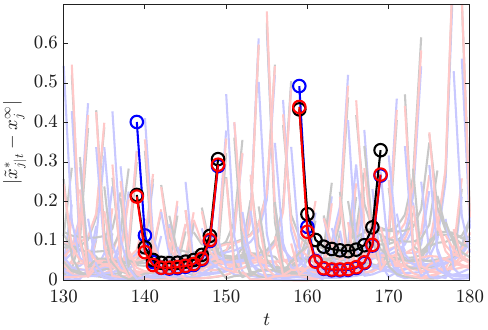}
	\caption{CSTR. Finite-horizon solutions $\tilde{x}^*_{j|t}$ compared to $x^\infty_j$, $j\in[t-N,t]$ for $t\in[130,180]$ using the filtering prior (blue), smoothing prior (black), and turnpike prior (red); highlighted are particular solutions obtained at $t=149$ and $t=169$.}
	\label{fig:MHE_CSTR_TP}
\end{figure}

\begin{figure}[t]
	\centering
	\includegraphics[]{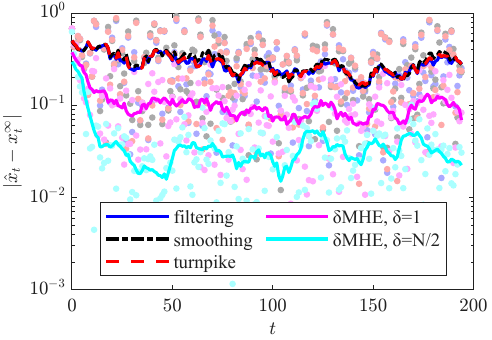}
	\caption{CSTR. Distance between state estimates using different MHE schemes and the IHE.
	Dots indicate values at time $t$, lines indicate the moving average over a sliding window of size $N+1$.}
	\label{fig:MHE_CSTR_E}
\end{figure}

From Figure~\ref{fig:MHE_CSTR_TP}, we can observe that all MHE problems (for all priors) exhibit the turnpike behavior with respect to the infinite-horizon solution, with clear approaching and leaving arcs, which is a strong indicator that Assumption~\ref{ass:turnpike_inf_prior} holds true.
We additionally compare standard MHE schemes (without delay) with $\delta$MHE~\eqref{eq:MHE_delay_x_prior} using the turnpike prior for $\delta=1$ and $\delta=N/2$.
From Figure~\ref{fig:MHE_CSTR_E}, we see that the standard MHE schemes yield very similar estimation results in terms of the difference to the IHE (for all priors), that $\delta$MHE with $\delta=1$ provides estimates that are much closer to the IHE, and that $\delta$MHE with $\delta=N/2$ ({which coincides with the turnpike prior as defined in~\eqref{eq:prior_turnpike}}) converges to a (small) neighborhood around the IHE, which nicely illustrates Proposition~\ref{prop:prior}.

\begin{figure}[t]
	\centering
	\includegraphics[]{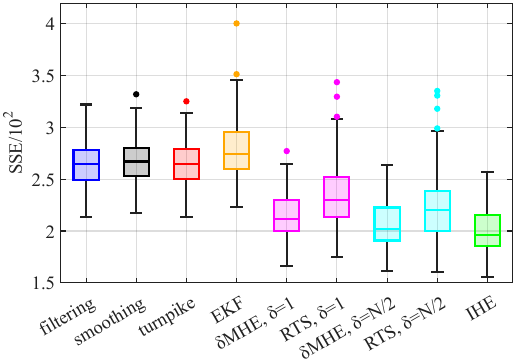}
	\caption{CSTR. Boxplot of the SSE for MHE using the filtering prior, the smoothing prior, and the turnpike prior, $\delta$MHE using the turnpike prior for $\delta=1$ and $\delta=N/2=5$, and the IHE over 100 different simulations with random disturbances and noise and randomly selected initial priors. {For comparison reasons, we also include the corresponding estimation results for the EKF and the RTS smoother using a fixed lag of $\delta=1$ and $\delta=N/2$.}}
	\label{fig:MHE_CSTR_SSE}
\end{figure}

We consider 100 different simulations with random disturbances and randomly selected initial priors $\bar{x}_0$ that are sampled from a uniform distribution over the interval centered at $x_0$ with a relative deviation of 25\% for each state.
Figure~\ref{fig:MHE_CSTR_SSE} shows that the standard MHE schemes are again very similar in terms of their SSE (for all priors), significantly outperformed by $\delta$MHE for $\delta=1$ (which yields a reduction in the SSE by 20\%). Moreover, we find that the SSE of $\delta$MHE with $\delta=N/2$ is very close to that of the IHE. {Finally, we observe that the optimization-based approaches outperform the EKF and corresponding fixed-lag RTS smoothers, which is due to the fact that these methods suffer from overshoots in some simulations that violate the state~constraints.}

Overall, this example nicely illustrates the developed theory. In particular, it shows that MHE problems with prior weighting exhibit the turnpike behavior with respect to the IHE (Assumption~\ref{ass:turnpike_inf_prior}), with a potentially strong leaving arc, cf.~Figure~\ref{fig:MHE_CSTR_TP}.
This motivates a delay in the estimation scheme in order to reduce the influence of the leaving arc.
Surprisingly, already a one-step delay is sufficient to significantly reduce the influence of the leaving arc such that $\delta$MHE tracks the performance and accuracy of the IHE with small error, cf.~Figures~\ref{fig:MHE_CSTR_E} and \ref{fig:MHE_CSTR_SSE}.

\subsubsection{Quadrotor}\label{sec:ex:MHE_quadrotor}

We consider a quadrotor with flexible blades and adapt the dynamical model from \cite{Kai2017}.
In the following, by $\mathcal{I}$ we denote the stationary inertial system with its vertical component pointing into the Earth, where $z=[z_1,z_2,z_3]^\top$ and $s=[s_1,s_2,s_3]^\top$ represent the position and velocity of the quadrotor, respectively.
By $\mathcal{B}$ we refer to the body-fixed frame attached to the quadrotor, with the third component pointing in the opposite direction of thrust generation.
The attitude of $\mathcal{B}$ with respect to $\mathcal{I}$ is described by a rotation matrix $R$ involving the roll, pitch, and yaw angle of the quadrotor, which we denote by $\xi=[\phi,\theta,\psi]^\top$.
The angular velocity of the quadrotor in $\mathcal{B}$ with respect to $\mathcal{I}$ is represented by $\Omega = [\Omega_1,\Omega_2,\Omega_3]^\top$.
The dynamics of the quadrotor can be described as
\begingroup
\renewcommand*{\arraystretch}{1.2}
\begin{align*}
	\begin{matrix*}[l]
		\dot{z}=s, &m\dot{s}=mge_3-TR(\xi)e_3-R(\xi)B\Omega,\\
		\dot{\xi}=\Gamma(\xi)\Omega, &J\dot{\Omega} = -\Omega^{\times}J\Omega+\tau-D\Omega,
	\end{matrix*}
\end{align*}
\endgroup
where $e_3 = [0\ 0\ 1]^\top$ and $(\cdot)^{\times}$ refers to the skew symmetric matrix associated with the cross product such that $u^{\times}v = u\times v$ for any $u,v\in\mathbb{R}^3$.
The thrust $T\in\mathbb{R}$ and the torque $\tau\in\mathbb{R}^3$ are generated by the velocities $\omega_i$, $i=\{1,2,3,4\}$ of the rotors via
\begin{equation*}
	\text{\small$\begin{bmatrix} T \\ \tau \end{bmatrix} = 
		\begin{bmatrix}
			c_T & c_T & c_T & c_T\\
			0 & -lc_T & 0 & lc_T\\
			lc_T & 0 & -lc_T & 0\\
			-c_Q & c_Q & -c_Q & c_Q
		\end{bmatrix}$}
	\text{\footnotesize$\begin{bmatrix} \omega_1^2 \\ \omega_2^2 \\ \omega_3^2 \\ \omega_4^2 \end{bmatrix}$},
\end{equation*}
and the matrix $\Gamma$ is defined as
\begin{equation*}
	\Gamma(\xi) = 
	\begin{bmatrix}
		1 & \sin\phi\tan\theta & \cos\phi\tan\theta\\
		0 & \cos\phi & -\sin\phi\\
		0 & \sin\phi\sec\theta & \cos\phi\sec\theta
	\end{bmatrix},
\end{equation*}
see, e.g.,~\cite{Kai2017} for modeling details and~\cite{Schiller2025_code} for implementation aspects, including parameter values and simulation settings.

We consider the state $x=[z^\top, \xi^\top, s^\top, \Omega^\top]^\top\in\mathbb{R}^{12}$ and input $u = [\omega_1, \omega_2, \omega_3, \omega_4]^\top\in\mathbb{R}^4$ and assume that only position $z$ and orientation $\xi$ are measured.
We derive the discrete-time model~\eqref{eq:sys} via the Euler method, yielding $f$ as in \eqref{eq:dynamics} and $h(x,u) = \left[ I_6, 0_{6\times6} \right]x$, while considering additional process disturbances $w\in\mathbb{R}^{12}$ and measurement noise $v\in\mathbb{R}^{6}$.

\begin{figure}[t]
	\centering
	\includegraphics[scale=0.73]{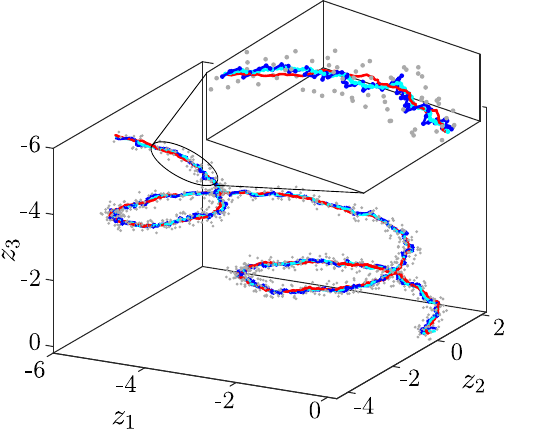}
	\caption{Quadrotor. Exemplary 3D trajectory for one specific disturbance realization; comparison of the true trajectory (red), measurements (gray dots), MHE with filtering prior (blue), and $\delta$MHE for $\delta=N/2=15$ (cyan).}
	\label{fig:MHE_quadrotor_X}
\end{figure}

In the simulations, $w$ and $v$ are uniformly distributed random variables sampled from bounded sets.
We consider the simulation time $T_\mathrm{s}=1000$ and a given open-loop trajectory $u_{0:T_\mathrm{s}}$, which moves the quadrotor spirally upwards, see Figure~\ref{fig:MHE_quadrotor_X} for an exemplary trajectory under a specific disturbance realization.
To estimate the unknown state $x_t$, we consider the cost function~\eqref{eq:J_p} with the horizon length $N=30$ and quadratic costs~\eqref{eq:stage_cost} and \eqref{eq:term_cost}, where we select $Q = \mathrm{blkdiag}(10^2I_3,10^4I_3,10^3I_3,10^5I_3)$ and  $R=G=\mathrm{blkdiag}(10^1I_3,10^2I_3)$.
Moreover, we use the quadratic prior weighting \eqref{eq:prior_cost}, where $W_t$ is initialized with $W_0=10I_{12}$ and updated using the EKF covariance formulas for all $t\in\mathbb{N}$.
We consider the case where no additional information about the domains of the states and disturbances is available and use $\mathcal{X}{\,=\,}\mathbb{R}^{12}$, $\mathcal{W}{\,=\,}\mathbb{R}^{12}$, $\mathcal{V}{\,=\,}\mathbb{R}^{6}$ in~\eqref{eq:MHE}.
In the following, we examine 100 different simulations with random disturbances and initial priors $\bar{x}_0$ sampled from a uniform distribution over $\mathcal{X}_0=\{x:|z_i|\leq 1, |\xi_i|\leq \pi/16, i=\{1,2,3\}, v=0, \Omega=0\}$.
We compare standard MHE with filtering prior~\eqref{eq:prior_filtering}, smoothing prior~\eqref{eq:prior_smoothing}, and turnpike prior~\eqref{eq:prior_turnpike}, $\delta$MHE~\eqref{eq:MHE_delay_x_prior} with turnpike-based prior weighting and $\delta=1$, $\delta=3$, $\delta=N/2$, and the IHE (which we approximate by solving the clairvoyant FIE problem using all available data $d_{0:T_{\mathrm{s}}}$).

\begin{figure}[t]
	\centering
	\includegraphics[]{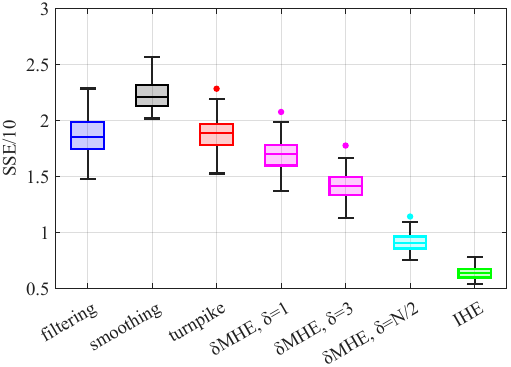}
	\caption{Quadrotor. Boxplot of the SSE for MHE using the filtering prior, the smoothing prior, and the turnpike prior, $\delta$MHE using the turnpike prior for $\delta=1$, $\delta=3$, and $\delta=N/2=15$, and the IHE over 100 different simulations with random disturbance/noise and randomly selected initial priors.}
	\label{fig:MHE_quadrotor_SSE}
\end{figure}

From Figure~\ref{fig:MHE_quadrotor_SSE}, we observe that MHE with filtering and turnpike prior performs quite similarly. The fact that MHE with smoothing prior is slightly worse can be attributed to the fact that the movement of the quadrotor is rather slow compared to the sampling time, while the horizon length $N=30$ is also rather small. In such setting, MHE with filtering or turnpike prior proves beneficial, as this essentially considers measurements from a larger estimation window, cf.~the discussion below~\eqref{eq:prior_smoothing}.
For $\delta$MHE, we can observe a reduction of the SSE of approximately 10\% for $\delta=1$, of 25\% for $\delta=3$, and of by 50\% for $\delta=N/2$, which is also close the SSE of the IHE.

To conclude, this example illustrates that the developed theory is also applicable to more complex and realistic systems.

\section{Conclusion}
In this paper, we developed novel accuracy and performance guarantees for optimal state estimation of general nonlinear systems, with a strong focus on MHE for online state estimation.
Our results rely on a turnpike property of the finite-horizon optimal state estimation problem with respect to the solution of the omniscient (acausal) infinite-horizon problem involving all past and future data.
This naturally causes MHE problems to exhibit a leaving arc that can have a potentially strong negative impact on estimation accuracy.
To counteract the leaving arc, we used an artificial delay in the MHE scheme, and we showed that the resulting performance is approximately optimal with respect to the infinite-horizon solution, with error terms that can be made arbitrarily small by an appropriate choice of the delay.
We proposed a novel turnpike prior for MHE formulations with prior weighting, proven to be a valid alternative to the classical options (such as the filtering or smoothing prior) with superior theoretical properties.

In our simulations, we found that MHE with the proposed turnpike prior performs comparably well to MHE with filtering or smoothing priors, while the delay resulted in a significant improvement of the estimation results.
In particular, considering practical examples from the literature, we observed the turnpike phenomenon and found that a delay of already one to three steps improved the overall estimation error by \mbox{20\,--25\%} compared to standard MHE (without delay).
For batch estimation, our proposed techniques have proven to be a useful alternative to classical filtering and smoothing methods.

{An interesting topic for future work is the application of the delayed MHE scheme in state feedback control---specifically, whether a delay in the closed loop is acceptable for obtaining significantly better estimation results. This could also motivate the development of alternative methods to reduce the influence of the leaving arc without using a delay, possibly based on available external predictions and forecasts.}


\begin{IEEEbiography}
	[{\includegraphics[width=1in,height=1.25in,clip,keepaspectratio]{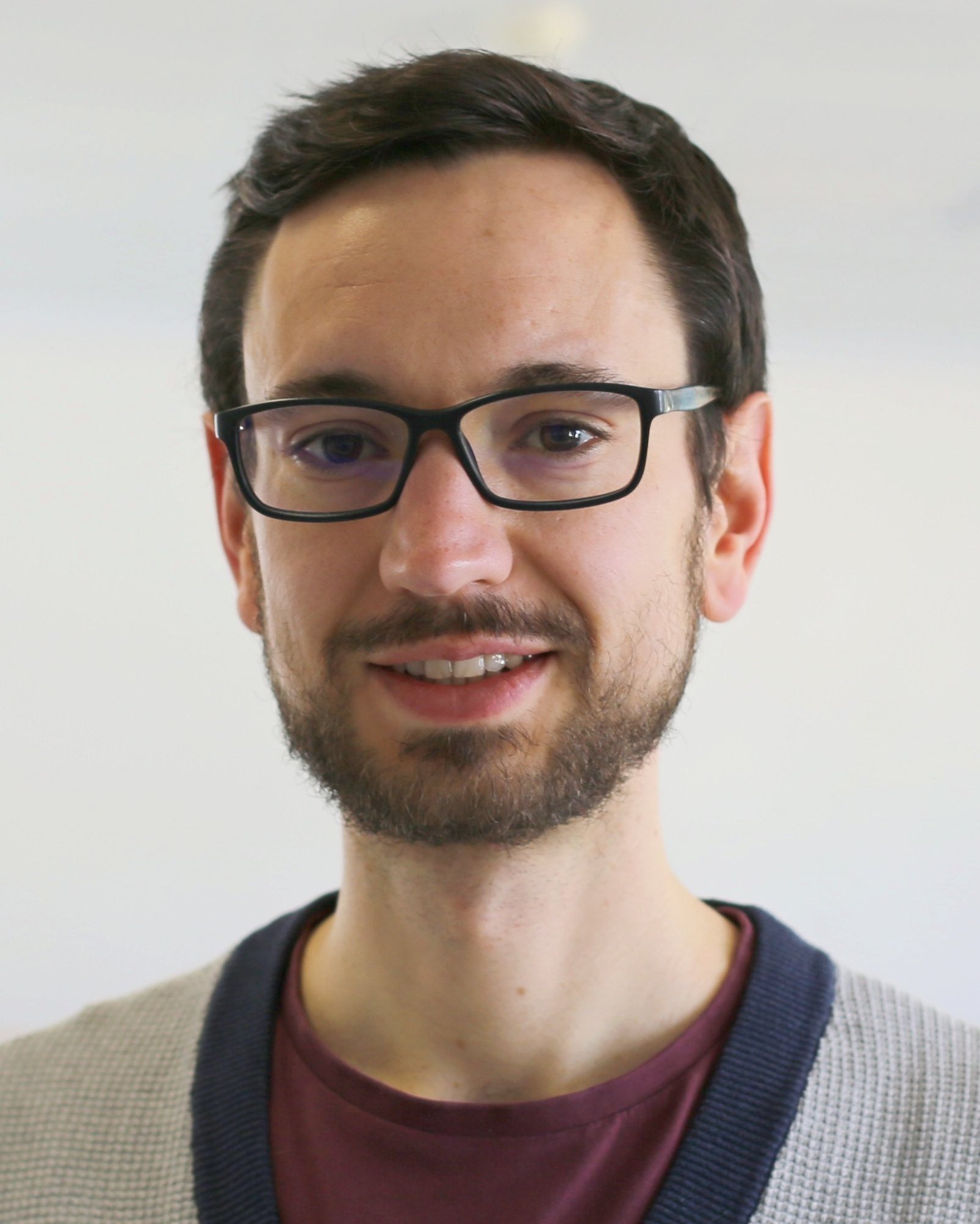}}]{Julian D. Schiller}
	received his Master degree in Mechatronics in 2019 and his Ph.D. degree in Electrical Engineering in 2025, both from the Leibniz University Hannover, Germany. He is currently a Postdoctoral Researcher at the Institute of Automatic Control, Leibniz University Hannover, Germany. His research interests are in the area of optimization-based state estimation and control of nonlinear systems and their interface with deep learning.
\end{IEEEbiography}

\begin{IEEEbiography}
	[{\includegraphics[width=1in,height=1.25in,clip,keepaspectratio]{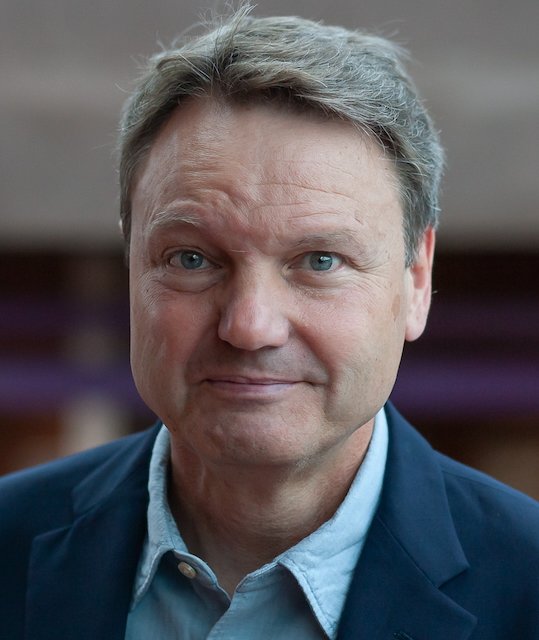}}]{Lars Grüne}
	has been Professor for Applied Mathematics at the University of Bayreuth, Germany, since 2002. He received his Diploma and Ph.D. in Mathematics in 1994 and 1996, respectively, from the University of Augsburg and his habilitation from the J.W. Goethe University in Frankfurt/M in 2001. He held or holds visiting positions at the Universities of Rome `Sapienza' (Italy), Padova (Italy), Melbourne (Australia), Paris IX - Dauphine (France), Newcastle (Australia) and IIT Bombay (India).
	
	Prof. Grüne was General Chair of the 25th International Symposium on Mathematical Theory on Networks and Systems (MTNS 2022), he is Editor-in-Chief of the journal Mathematics of Control, Signals and Systems (MCSS) and is or was Associate Editor of various other journals, including the Journal of Optimization Theory and Applications (JOTA), Mathematical Control and Related Fields (MCRF) and the IEEE Control Systems Letters (CSS-L). His research interests lie in the area of mathematical systems and control theory with a focus on numerical and optimization-based methods for nonlinear systems.
\end{IEEEbiography}

\begin{IEEEbiography}
	[{\includegraphics[width=1in,height=1.25in,clip,keepaspectratio]{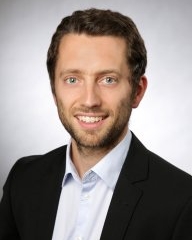}}]{Matthias A. M\"uller}
	received a Diploma degree in engineering cybernetics from the University of Stuttgart, Germany, an M.Sc. in electrical and computer engineering from the University of Illinois at Urbana-Champaign, US (both in 2009), and a Ph.D. from the University of Stuttgart in 2014. Since 2019, he is Director of the Institute of Automatic Control and Full Professor at the Leibniz University Hannover, Germany. 
	
	His research interests include nonlinear control and estimation, model predictive control, and data- and learning-based control, with application in different fields including biomedical engineering and robotics. He has received various awards for his work, including the 2015 European Systems \& Control PhD Thesis Award, the inaugural Brockett-Willems Outstanding Paper Award for the best paper published in Systems \& Control Letters in the period 2014-2018, an ERC starting grant in 2020, the IEEE CSS George S. Axelby Outstanding Paper Award 2022, and the Journal of Process Control Paper Award 2023. He serves/d as an associate editor for Automatica and as an editor of the International Journal of Robust and Nonlinear Control.
	%
\end{IEEEbiography}

\end{document}